\documentclass{article}
 \usepackage{color}
\usepackage{amsmath}
\usepackage{amsthm}
\usepackage{amsfonts}
\usepackage{amssymb}
\usepackage{esint}
\usepackage{eucal}
\usepackage{mathrsfs}
\usepackage{graphicx,graphics}  
\numberwithin{equation}{section}
\usepackage{graphicx,graphics}
\usepackage{pdfsync}
\usepackage{multirow} 
\usepackage{endnotes}
\usepackage{epstopdf}
\usepackage{color}
\usepackage{enumerate}
\usepackage[colorlinks=true, linkcolor=blue]{hyperref} % For coloured hyperrefs

\usepackage{datetime}
\def \dis {\displaystyle}

\def \into {\int_\Omega}
\def \confai {-\kern -.5em\rightharpoonup}
\def \cqfd {\hfill$\Box$}

\def \al {\alpha}
\def \be {\beta}

\def \Ga {\Gamma}

\def \om {\omega}
\def \Om {\Omega}
\def \la {\lambda}

\def \ph {\varphi}

\def \NN {\mathbb N}

\def \RR {\mathbb R}

\def \D {\mathscr{D}}
\def \F {\mathscr{F}}

\def \L {\mathscr{L}}
\def \M {\mathscr{M}}

\def \beq {\begin{equation}}
\def \eeq {\end{equation}}
\def \ba {\begin{array}}
\def \ea {\end{array}}
\def \bs {\bigskip}
\def \ms {\medskip}

\def \ecart {\noalign{\medskip}}

\newtheorem{theorem}{Theorem}[section]
\newtheorem{lemma}[theorem]{Lemma}
\newtheorem{corollary}[theorem]{Corollary}
\newtheorem{definition}[theorem]{Definition}
\newtheorem{remark}[theorem]{Remark}

\DeclareMathOperator{\Dive}{Div}

\DeclareMathOperator{\dist}{dist}
\DeclareMathOperator*{\essinf}{ess-inf}

\title{Homogenization of equi-coercive nonlinear energies defined on vector-valued functions, with non-uniformly bounded coefficients}
\author{
M.~Briane\footnote{IRMAR \& INSA Rennes, mbriane@insa-rennes.fr},
J.~Casado-D\'iaz\footnote{Dpto. de Ecuaciones Diferenciales y An\'alisis Num\'erico, Universidad de Sevilla, jcasadod@us.es},
M.~Luna-Laynez\footnote{Dpto. de Ecuaciones Diferenciales y An\'alisis Num\'erico, Universidad de Sevilla, mllaynez@us.es},
A.~Pallares-Mart\'in\footnote{Dpto. de Ecuaciones Diferenciales y An\'alisis Num\'erico, Universidad de Sevilla, ajpallares@us.es}
}

\begin{document}

\maketitle
\par\bigskip\noindent

\begin{abstract}
The present paper deals with the asymptotic behavior of equi-coercive sequences $\{\F_n\}$ of nonlinear functionals defined over vector-valued functions in $W^{1,p}_0(\Om)^M$, where $p>1$, $M\geq 1$, and $\Om$ is a bounded open set of $\RR^N$, $N\geq 2$. The strongly local energy density $F_n(\cdot,Du)$ of the functional $\F_n$ satisfies a Lipschitz condition with respect to the second variable, which is controlled by a positive sequence $\{a_n\}$ which is only bounded in some suitable space $L^r(\Om)$. We prove that the sequence $\{\F_n\}$ $\Gamma$-converges for the strong topology of $L^p(\Om)^M$ to a functional $\F$ which has a strongly local density $F(\cdot,Du)$ for sufficiently regular functions~$u$. This compactness result extends former results on the topic, which are based either on maximum principle arguments in the nonlinear scalar case, or adapted div-curl lemmas in the linear case. Here, the vectorial character and the nonlinearity of the problem need   a new approach based on a careful analysis of the asymptotic minimizers associated with the functional $\F_n$. The relevance of the conditions which are imposed to the energy density $F_n(\cdot,Du)$, is illustrated by several examples including some classical hyper-elastic energies.
\end{abstract}
\section{Introduction}
In this paper we study the asymptotic behavior of the sequence of nonlinear functionals, including some hyper-elastic energies (see the examples of Section~\ref{hyperelas}), defined on vector-valued functions by
\beq\label{Fni}
\F_n(v):=\int_\Omega F_n(x,Dv)\,dx\quad\mbox{for }v\in W^{1,p}_0(\Om)^M,\quad\mbox{with }p\in(1,\infty),\ M\geq 1,
\eeq
in a bounded open set $\Om$ of $\RR^N$, $N\geq 2$.
{ 
The sequence $\F_n$ is assumed to be equi-coercive. Moreover, the associated density~$F_n(\cdot,\xi)$ satisfies some Lipschitz condition with respect to $\xi\in\RR^{M\times N}$, and its coefficients are not uniformly bounded in~$\Om$.
\par
The linear scalar case, {\em i.e.} when $F_n(\cdot,\xi)$ is quadratic with respect to $\xi\in\RR^N$ ($M=1$), with uniformly bounded coefficients was widely investigated in the seventies through G-convergence by Spagnolo \cite{Spa}, extended by Murat and Tartar with H-convergence~\cite{Mur,Tar}, and alternatively through $\Gamma$-convergence by De Giorgi \cite{DeG1,DGFr} (see also \cite{Dal,Bra}). The linear elasticity case was probably first derived by Duvaut (unavailable reference), and can be found in \cite{San,Fra}. In the nonlinear scalar case the first compactness results are due to Carbone, Sbordone \cite{CaSb} and Buttazzo, Dal Maso \cite{BuDa} by a $\Gamma$-convergence approach assuming the $L^1$-equi-integrability of the coefficients. More recently, these results were extended in \cite{BBC,BrCa1,BrCa3} relaxing the $L^1$-boundedness of the coefficients but assuming that $p>N\!-\!1$ if $N\geq 3$, showing then the uniform convergence of the minimizers thanks to the maximum principle. In all these works the scalar framework combined with the condition $p>N\!-\!1$ if $N\geq 3$ and the equi-coercivity of the functionals, induce in terms of the $\Gamma$-convergence for the strong topology of $L^p(\Om)$, a limit energy $\F$ of the same nature satisfying
\beq\label{Fi}
\F(v):=\int_\Omega F(x,Dv)\,d\nu\quad\mbox{for }v\in W,
\eeq
where $C^1_c(\Om)^M\subset W$ is some suitable subspace of $W^{1,p}_0(\Om)^M$, and $\nu$ is some Radon measure on $\Om$.
Removing the $L^1$-equi-integrability of the coefficients in the three-dimensional linear scalar case (note that $p=N\!-\!1=2$ in this case),
Fenchenko and Khruslov \cite{FeKh} (see also~\cite{Khr}) were, up to our knowledge,}
the first to obtain a violation of the compactness result due to the appearance of local and nonlocal terms in the limit energy $\F$. This seminal work was also revisited by Bellieud and Bouchitt\'e \cite{BeBo}. Actually, the local and nonlocal terms in addition to the classical strongly local term come from the Beurling-Deny \cite{BeDe} representation formula of a Dirichlet form, and arise naturally in the homogenization process as shown by Mosco~\cite{Mos}. The complete picture of the attainable energies was obtained by Camar-Eddine and Seppecher~\cite{CESe1} in the linear scalar case. The elasticity case  is much more intricate even in the linear framework, since the loss of uniform boundedness of the elastic coefficients may induce the appearance of second gradient terms as Seppecher and Pideri proved in~\cite{PiSe}. The situation is dramatically different from the scalar case, since the Beurling-Deny formula does not hold in the vector-valued case. In fact, Camar-Eddine and Seppecher~\cite{CESe2} proved that any lower semi-continuous quadratic functional vanishing on the rigid displacements, can be attained. Compactness results were obtained in the linear elasticity case using some (strong) equi-integrability of the coefficients in \cite{BrCa4}, and using various extensions of the classical Murat-Tartar~\cite{Mur} div-curl result in \cite{BrCE1, BCM, BrCa5, APM} (which were themselves initiated in the former works \cite{Bri,BrCa1} of the two first authors).
\par
In our context the vectorial character of the problem and its nonlinearity prevent us from using the uniform convergence of \cite{BrCa3} and the div-curl lemma of \cite{BrCa5}, which are (up to our knowledge) the more recent general compactness results on the topic.
We assume that the nonnegative energy density $F_n(\cdot,\xi)$ of the functional~\eqref{Fni} attains its minimum at $\xi=0$, and satisfies the following Lipschitz condition with respect to $\xi\in\RR^{M\times N}$:
\begin{equation*}%\label{cLip}
 \left\{\ba{ll}
    \big|F_n(x,\xi)-F_n(x,\eta)\big|\leq\big(h_n(x)+F_n(x,\xi)+F_n(x,\eta)+|\xi|^p+|\eta|^p\big)^{\frac{p-1}{p}}a_n(x)^{\frac{1}{p}}\,|\xi-\eta|
    \\ \ecart
    \forall\,\xi,\eta \in \mathbb{R}^{M\times N},\mbox{ a.e. }x\in \Omega,
  \ea\right.
\end{equation*}
which is controlled by a positive function $a_n(\cdot)$ (see the whole set of conditions \eqref{hip1} to \eqref{hipo3} below). The sequence $\{a_n\}$ is assumed to be bounded in $L^r(\Om)$ for some $r>(N\!-\!1)/p$ if $1<p\leq N\!-\!1$, and bounded in $L^1(\Om)$ if $p>N\!-\!1$. Note that for $p>N\!-\!1$ our condition is better than the $L^1$-equi-integrability used in the scalar case of~\cite{CaSb,BuDa}, but not for $1<p\leq N\!-\!1$. Under these assumptions we prove (see Theorem~\ref{thm2}) that the sequence $\{\F_n\}$ of~\eqref{Fni} $\Gamma$-converges for the strong topology of $L^p(\Om)^M$ (see Definition~\ref{defGc}) to a functional of type \eqref{Fi} with
\begin{equation*}%\label{Wm}
W\subset\left\{\ba{ll}
W^{1,{pr\over r-1}}(\Om)^M, & \mbox{if }1<p\leq N\!-\!1
\\ \ecart
C^1(\overline{\Om})^M, & \mbox{if }p>N\!-\!1,
\ea\right.
\quad\mbox{and}\quad
\nu=\left\{\ba{ll}
\dis \mbox{Lebesgue measure}, & \mbox{if }1<p\leq N\!-\!1
\\ \ecart
\dis {\M(\Om)\,*}-\lim_{n\to\infty}a_n, & \mbox{if }p>N\!-\!1.
\ea\right.
 \end{equation*}
 Various types of boundary conditions can be taken into account in this $\Gamma$-convergence approach.
 \par
 A preliminary result (see Theorem~\ref{thm}) allows us to prove that the sequence of energy density $\{F_n(\cdot,Du_n)\}$ converges in the sense of Radon measures to some strongly local energy density $F(\cdot,Du)$, when $u_n$ is an asymptotic minimizer for $\F_n$  of limit $u$ (see definition~\eqref{reco1bis}). The proof of this new compactness result is based on an extension (see Lemma~\ref{lemafund}) of the fundamental estimate for recovery sequences in $\Gamma$-convergence (see, {\em e.g.},~\cite{Dal}, Chapters~18,~19), which provides a bound (see \eqref{acopmunu}) satisfied by the weak-$*$ limit of $\{F_n(\cdot,Du_n)\}$ with respect to the weak-$*$ limit of any sequence $\{F_n(\cdot,Dv_n)\}$ such that the sequence $\{v_n-u_n\}$ converges weakly to $0$ in $W^{1,p}_0(\Om)^M$.
Rather than using fixed smooth cut-off functions as in the classical fundamental estimate, here we need to consider sequences of radial cut-off functions $\ph_n$ whose gradient has support in $n$-dependent sets on which $u_n-u$ satisfies some uniform estimate with respect to the radial coordinate (see Lemma~\ref{lem1} and its proof). This allows us to control the zero-order term $\nabla\ph_n(u_n-u)$, when we put the trial function $\ph_n(u_n-u)$ in the functional $\F_n$ of~{\eqref{Fni}. The uniform estimate is a consequence of the Sobolev compact embedding for the $(N\!-\!1)$-dimensional sphere, and explains the role of the exponent $r>(N\!-\!1)/p$ if $1<p\leq N\!-\!1$. A similar argument was used in the linear case~\cite{BrCa5} to obtain a new div-curl lemma which is the key-ingredient for the compactness of quadratic elasticity functionals of type~\eqref{Fni}.
\subsection*{Notations}
\begin{itemize}
{ 
\item $\RR_s^{N\times N}$ denotes the set of the symmetric matrices in $\RR^{N\times N}$.
\item For any $\xi\in\RR^{N\times N}, \xi^T$ is the transposed matrix of $\xi$, and $\xi^s:={1\over 2}(\xi+\xi^T)$ is the symmetrized matrix of $\xi$.
\item $I_N$ denotes the unit matrix of $\RR^{N\times N}$.
\item $\cdot$ denotes the scalar product in $\RR^N$, and $:$ denotes the scalar product in $\RR^{M\times N}$ defined by
\par
\centerline{$\xi:\eta:={\rm tr}\,(\xi^T\eta)\quad$for $\xi,\eta\in\RR^{M\times N}$,}
where tr is the trace.
}
\item $|\cdot|$ denotes both the euclidian norm in $\RR^N$, and the Frobenius norm in $\RR^{M\times N}$, {\em i.e.}
\par
\centerline{$|\xi|:=\big({\rm tr}\,(\xi^T\xi)\big)^{1\over 2}\quad$for $\xi\in\RR^{M\times N}$.}
\item {  For a bounded open set $\om\subset \RR^N$, $\M(\om)$ denotes the space of the Radon measures on $\om$ with bounded total variation. It agrees with the dual space of $C^0_0(\om)$, namely the space of the continuous functions in $\bar \om$ which vanish on $\partial\om$. Moreover, $\M(\bar \om)$ denotes the space of the Radon measures on $\bar \om$. It agrees with the dual space of $C^0(\bar\om)$.}
\item For any measures $\zeta, \mu\in\M(\omega)$, with $\om\subset\RR^N$, open, bounded, we define $\zeta^\mu\in L^1_\mu(\Omega)$ as the derivative of $\zeta$ with respect to $\mu$. When $\mu$ is the Lebesgue measure, we write  $\zeta^L$.
\item $C$ is a positive constant which may vary from line to line.
\item $O_n$ is a real sequence which tends to zero as $n$ tends to infinity. It can vary from line to line.
\end{itemize}
Recall the definition of the De Giorgi $\Gamma$-convergence (see, {\em e.g.}, \cite{Dal,Bra} for further details).
\begin{definition}\label{defGc}
Let $V$ be a metric space, and let $\F_n,\F:V\to [0,\infty]$, $n\in\NN$, be functionals defined on $V$. The sequence $\{\F_n\}$ is said to $\Ga$-converge to $\F$ for the topology of $V$ in a set $W\subset V$ and we write
$$\F_n\stackrel{\Gamma}\rightharpoonup \F\ \hbox{ in }W,$$ 
if
\begin{itemize}
\item[-] the $\Ga$-liminf inequality holds
\begin{equation*}
\label{Gli}
{\forall\,v\in W,}\  \ \forall\,v_n\rightarrow v\;\;\mbox{ in }V,\quad \F(v)\leq\liminf_{n\to\infty}\F_n(v_n),
\end{equation*}
\item[-] the $\Ga$-limsup inequality holds
\begin{equation*}
\label{Gls}
{\forall v\,\in W,}\ \ \exists\,\overline{v}_n \rightarrow  v\;\;\mbox{ in }V,\quad \F(v)=\lim_{n\to\infty}\F_n(\overline{v}_n).
\end{equation*}
\end{itemize}
Any sequence $\overline{v}_n$ satisfying~\eqref{Gls} is called a {\em recovery sequence} for $\F_n$ of limit $v$.
\end{definition}

\section{Statement of the results and examples}
\subsection{The main results}
Consider a bounded open set $\Omega\subset \mathbb{R}^N$ with $N\geq 2$, $M$ a positive integer, a sequence of nonnegative Carath\'{e}odory functions $F_n:\Omega\times \mathbb{R}^{M\times N}\to [0,\infty)$, and $p>1$ with the following properties:
\begin{itemize}
\item There exist two constants $\alpha>0$ and $\be\in\RR$ such that
\begin{equation}\label{hip1}
  \int_\Omega F_n(x,Du)\,dx\geq \alpha \int_\Omega |Du|^p\,dx+\be,\quad\forall\, u\in W^{1,p}_0(\Omega)^M,
\end{equation}
and
\begin{equation}\label{acoL1}
F_n(\cdot,0)=0\ \hbox{ a.e. in }\Om.
\end{equation}
\item There exist two sequences of  measurable functions $h_n, a_n\geq0$, and a constant $\gamma>0$ such that

\begin{equation}\label{acohn}
  h_n \mbox{ is bounded in }L^1(\Omega),
\end{equation}

\begin{equation}\label{acoan}
  a_n \mbox{ is bounded in }L^r(\Omega) \mbox{ with }
    \left\{\ba{ll}
        \dis r > \frac{N\!-\!1}{p}, & \mbox{if }1<p\leq N\!-\!1
        \\ \ecart
       r=1, & \mbox{if } p > N\!-\!1,
    \ea\right.
\end{equation}

\begin{equation}\label{hipo2}
  \left\{\ba{ll}
    \big|F_n(x,\xi)-F_n(x,\eta)\big|\leq\big(h_n(x)+F_n(x,\xi)+F_n(x,\eta)+|\xi|^p+|\eta|^p\big)^{\frac{p-1}{p}}a_n(x)^{\frac{1}{p}}\,|\xi-\eta|
    \\ \ecart
    \forall\,\xi,\eta \in \mathbb{R}^{M\times N},\mbox{ a.e. }x\in \Omega,
  \ea\right.
\end{equation}
and
\begin{equation}\label{hipo3}
  F_n(x,\lambda\xi)\leq h_n(x)+\gamma\,F_n(x,\xi),
  \quad\forall\,\lambda\in [0,1],\ \forall\,\xi\in\mathbb{R}^{M\times N},\mbox{ a.e. }x\in\Omega.
\end{equation}
\end{itemize}

\begin{remark}
  From \eqref{hipo2} and Young's inequality, we get that
  \begin{align*}
    F_n(x,\xi) & \leq F_n(x,\eta)+\big(h_n(x)+F_n(x,\xi)+F_n(x,\eta)+|\xi|^p+|\eta|^p\big)^{\frac{p-1}{p}}a_n(x)^{\frac{1}{p}}\,|\xi-\eta|
    \\
    & \leq F_n(x,\eta)+\frac{p-1}{p}\big(h_n(x)+F_n(x,\xi)+F_n(x,\eta)+|\xi|^p+|\eta|^p\big)+\frac{1}{p}\,a_n(x)\,|\xi - \eta|^p,
  \end{align*}
  and then
  \begin{equation}\label{acotFnb}
    F_n(x,\xi)\leq (p-1)\,h_n(x)+(2p-1)\,F_n(x,\eta)+(p-1)\big(|\xi|^p+|\eta|^p\big)+a_n(x)\,|\xi-\eta|^p,\quad\forall\, \xi,\eta\in \mathbb{R}^{M\times N},\mbox{ a.e. }x\in\Omega.
  \end{equation}

  In particular, taking $\eta = 0$, we have
  \begin{equation}\label{acotFn}
    F_n(x,\xi)\leq (p-1)\,h_n(x)+\big(p-1+a_n(x)\big)\,|\xi|^p,\quad\forall\,\xi\in \mathbb{R}^{M\times N},\mbox{ a.e. } x\in\Omega,
  \end{equation}
\noindent
where the right-hand side is a bounded sequence in $L^1(\Omega)$.\newline
\end{remark}
  From now on, we assume that
  \begin{equation}\label{deliahn}
     a_n^r\stackrel{*} \rightharpoonup {\textsc{a}}\;\;\mbox{in }\M(\Omega)\quad\mbox{and}\quad
     h_n\stackrel{*}\rightharpoonup {h}\;\;\mbox{in }\M(\Omega).
  \end{equation}

The paper deals with the asymptotic behavior of the sequence of functionals
\beq\label{FnOm}
\F_n(v):=\int_\Omega F_n(x,Dv)\,dx\quad\mbox{for }v\in W^{1,p}(\Om)^M.
\eeq

First of all, we have the following result on the convergence of the energy density $F_n(\cdot,D u_n)$, where $u_n$ is an asymptotic minimizer associated with functional \eqref{FnOm}.
   \begin{theorem}\label{thm}
Let $F_n:\Omega\times \mathbb{R}^{M\times N}\to [0,\infty)$ be a sequence of Carath\'{e}odory functions satisfying \eqref{hip1} to \eqref{hipo3}. Then, there exist a function $F:\Omega\times\mathbb{R}^{M\times N}\to\mathbb{R}$ and a subsequence of $n$, still denoted by $n$, such that for any $\xi,\eta\in\RR^N$,
 \beq\label{proF0}
 \left\{\ba{ll}\dis F(\cdot,\xi)\ \hbox{ is Lebesgue measurable}, &\hbox{ if }1<p\leq N\!-\!1
 \\ \ecart
 \dis F(\cdot,\xi)\ \hbox{ is }{\textsc{a}}\hbox{-measurable}, &\dis\hbox{ if }p>N\!-\!1,\ea\right.
 \eeq

\begin{equation}\label{proF1}\ba{l}  \big|F(x,\xi)-F(x,\eta)\big|\leq
\\ \ecart
\dis \left\{\ba{ll}\dis C\big({h}^L\hskip-1pt+\hskip-1pt F(x,\xi)\hskip-1pt+\hskip-1pt F(x,\eta)\hskip-1pt+\hskip-1pt(1+({\textsc{a}}^L)^\frac{1}{r})(|\xi|^p+|\eta|^p)\big)^\frac{p-1}{p}\hskip-4pt({\textsc{a}}^L)^\frac{1}{pr}\hskip-1pt|\xi-\eta|\hbox{ a.e. in }\Om, &\hbox{if }1<p\leq\hskip-1pt N \hskip-2pt -\hskip-2pt 1
\\ \ecart
\dis C\big(1+{h} ^{\textsc{a}}+F(x,\xi)+F(x,\eta)+|\xi|^p+|\eta|^p\big)^\frac{p-1}{p}|\xi-\eta|\ \ {\textsc{a}}\hbox{-a.e. in }\Om, &\hbox{if }p>\hskip-1pt N \hskip-2pt -\hskip-2pt 1,\ea\right.
\ea\eeq
and
\beq\label{proF2} F(\cdot,0)=0\;\;\hbox{a.e. in }\Om.
\eeq

For any open set $\om\subset \Omega$, and any sequence $\{u_n\}$ in $W^{1,p}(\om)^M$ which converges weakly in $W^{1,p}(\om)^M$ to a function $u$ satisfying

\beq\label{regulu} u\in \left\{\ba{ll}\dis W^{1,\frac{pr}{r-1}}(\om)^M, & \hbox{ if }1< p\leq N\!-\!1
\\ \ecart
\dis C^1(\om)^M, &\hbox{ if }p>N\!-\!1,\ea\right.
\eeq
and such that
  \begin{equation}\label{reco1bis}
    \exists \lim_{n\to\infty}\int_\om F_n(x,Du_n)\,dx=\min\left\{\liminf_{n\to\infty} \int_\om F_n(x,Dw_n)\,dx:\ w_n-u_n\rightharpoonup 0\mbox{ in }W^{1,p}_0(\om)^M\right\}<\infty,
  \end{equation}
we have
\beq\label{concthp}
F_n(\cdot,Du_n)\stackrel{*}\rightharpoonup
\left\{\ba{ll}
F(\cdot,Du), & \hbox{if }1<p\leq N\!-\!1
\\ \ecart
F(\cdot,Du)\,{\textsc{a}}, & \mbox{if }p>N\!-\!1
\ea\right.
\quad\hbox{in }\M(\omega).
\eeq
\end{theorem}
\par\bs
From Theorem~\ref{thm} we may deduce the $\Gamma$-limit (see Definition~\ref{defGc}) of the sequence of functionals \eqref{FnOm} with various boundary conditions.
\begin{theorem}\label{thm2}
Let $F_n:\Omega\times \mathbb{R}^{M\times N}\to [0,\infty)$ be a sequence of Carath\'{e}odory functions satisfying \eqref{hip1} to \eqref{hipo3}. Let $\om$ be an open set such that $\om\subset\subset\Om$, and let $V$ be a subset of $W^{1,p}(\om)^M$ such that
\beq\label{V}
\forall\,u\in V,\ \forall\,v\in W^{1,p}_0(\om)^M,\quad u+v\in V.
\eeq

Define the functional $\F_n^V:V\to [0,\infty)$ by
\beq\label{FnV}
\F_n^V(v):=\int_{\omega} F_n(x,Dv)\,dx\quad\mbox{for }v\in V.
\eeq
Assume that the open set $\om$ satisfies
\beq\label{omega}
\left\{\ba{ll}
|\partial\om|=0, & \mbox{if }1<p\leq N\!-\!1
\\ \ecart
{\textsc{a}}(\partial\om)=0, & \mbox{if }p>N\!-\!1.
\ea\right.
\eeq
\noindent
Then, for the subsequence of $n$ (still denoted by $n$) obtained in Theorem~\ref{thm} we get {\beq\label{FV}\left\{\ba{ll}\dis \F_n^V\stackrel{\Gamma}\rightharpoonup \F^V:=\int_\om F(x,Dv)\,dx\  \mbox{ in } V\cap W^{1,{pr\over r-1}}( \om)^M, & \mbox{if }1<p\leq N\!-\!1\\ \ecart\dis
\F_n^V\stackrel{\Gamma}\rightharpoonup \F^V:=\int_\om F(x,Dv)\,dx\  \mbox{ in } V\cap C^1(\overline{\om})^M, & \mbox{if }p>N\!-\!1,
\ea\right.
\eeq
for the strong topology of $L^p(\om)^M$}, where $F$ is given by convergence~\eqref{concthp}.
\end{theorem}

\begin{remark}\label{remthm2}
The condition \eqref{omega} on the open set $\omega$ is not so restrictive. Indeed, for any family $(\om)_{i\in I}$ of open sets of $\Om$ with  two by two disjoint boundaries, at most a countable subfamily of $(\partial\om)_{i\in I}$ does not satisfy~\eqref{omega}.
\end{remark}
\subsection{Auxiliary lemmas}
The proof of Theorem~\ref{thm} is based on the following lemma which provides an estimate of the energy density for asymptotic minimizers. In our context it is equivalent to the fundamental estimate for recovery sequences (see Definition~\ref{defGc}) in $\Gamma$-convergence theory (see, {\em e.g.}, \cite{Dal}, Chapters~18,~19).
\begin{lemma}\label{lemafund}
Let $F_n:\Omega\times \mathbb{R}^{M\times N}\to [0,\infty)$ be a sequence of Carath\'{e}odory functions satisfying \eqref{hip1} to \eqref{hipo3}. Consider an open set $\om\subset \Omega$, and a sequence $\{u_n\}\subset W^{1,p}(\om)^M$ converging weakly in $W^{1,p}(\om)^M$ to a function $u$ satisfying \eqref{regulu}, and such that
\begin{equation*}%\label{comedu}
    F_n({\cdot},Du_n)\stackrel{*}\rightharpoonup \mu\quad\mbox{in }\M(\omega),
  \end{equation*}
\begin{equation*}%\label{devr}
|Du_n|^p\stackrel{*}\rightharpoonup \varrho\quad\mbox{in }\M(\omega).
\end{equation*}
Then, the measure $\varrho$ satisfies
\beq\label{estrho}
\varrho\leq\left\{\ba{lll}
C\big(|Du|^p+|Du|^p({\textsc{a}}^L)^{1\over r}+{h}+\mu+{\textsc{a}}^L\big) & \mbox{a.e. in }\om, & \mbox{if }1<p\leq N\!-\!1
\\ \ecart
C\big(|Du|^p {\textsc{a}}+{h}+\mu+{\textsc{a}}\big) & {\textsc{a}}\mbox{-a.e. in }\om, & \mbox{if }p>N\!-\!1.
\ea\right.
\eeq
\noindent
{Moreover if $u_n$ satisfies}
\begin{equation}\label{reco1}
    \exists \lim_{n\to\infty}\int_\om F_n(x,Du_n)\,dx=\min\left\{\liminf_{n\to\infty} \int_\om F_n(x,Dw_n)\,dx:\ w_n-u_n\rightharpoonup 0\mbox{ in }W^{1,p}_0(\omega)^M\right\},
  \end{equation}
then for any sequence $\{v_n\}\subset W^{1,p}(\omega)^M$ which converges weakly in $W^{1,p}(\omega)^M$ to a function
\[
v\in \left\{\ba{ll}\dis W^{1,\frac{pr}{r-1}}(\om)^M, &\hbox{ if }1< p\leq N\!-\!1\\ \ecart\dis C^1(\om)^M, &\hbox{ if }p>N\!-\!1,\ea\right.
\]
and such that
\begin{equation}\label{comedv}
F_n({\cdot},Dv_n)\stackrel{*}\rightharpoonup \nu\quad\mbox{in }\M(\om),
\end{equation}
  \begin{equation*}%\label{comedvb}
    |Dv_n|^p\stackrel{*}\rightharpoonup \varpi \quad\mbox{in }\M(\om),
  \end{equation*}
we have
\begin{equation}\label{acopmunu}
\mu\leq
\begin{cases}
\dis \nu+C\big({h}^L+\nu^L+\varpi^L+(1+({\textsc{a}}^L)^\frac{1}{r})|D(u-v)|^p \big)^\frac{p-1}{p}({\textsc{a}}^L)^\frac{1}{pr}|D(u-v)|\mbox{ a.e. in }\om, & \hbox{if }1<p\leq N\!-\!1
\\ \ecart
\dis \nu+C\big(1+{h} ^{\textsc{a}}+\nu ^{\textsc{a}}+\varpi ^{\textsc{a}}+|D(u-v)|^p\big)^\frac{p-1}{p}{\textsc{a}}\,|D(u-v)|\ {\textsc{a}}\mbox{-a.e. in }\om, & \hbox{if }p>N\!-\!1.
\end{cases}
\end{equation}
\end{lemma}

We can {improve} the statement of Lemma~\ref{lemafund} {if we add} a non-homogeneous Dirichlet boundary condition on {$\partial\om$}.
 \begin{lemma}\label{lemafund2}
Let $\om$ be an open set such that $\om\subset\subset\Om$, and let $u$ be a function satisfying
\beq\label{regulu1}
u\in \left\{\ba{ll}\dis W^{1,\frac{pr}{r-1}}(\Om)^M, & \hbox{ if }1< p\leq N\!-\!1
\\ \ecart
\dis C^1(\overline{\Om})^M, &\hbox{ if }p>N\!-\!1.
\ea\right.
\eeq
Let $\{u_n\}$ and $\{v_n\}$ be two sequences in $W^{1,p}(\om)^M$, {  such that $u_n$ satisfies condition \eqref{reco1} and}
\begin{equation*}
{u_n-u,\ v_n-u\in W^{1,p}_0(\omega)^M,}
\end{equation*}
\begin{equation}\label{comedu2}
F_n({\cdot},Du_n)\stackrel{*}\rightharpoonup \mu\quad\mbox{and}\quad F_n({\cdot},Dv_n)\stackrel{*}\rightharpoonup \nu
\quad\mbox{in }\M(\overline{\om}),
\end{equation}
\begin{equation}\label{devr2}
|Du_n|^p\stackrel{*}\rightharpoonup \varrho\quad\mbox{and}\quad
|Dv_n|^p\stackrel{*}\rightharpoonup \varpi\quad\mbox{in }\M(\overline{\om}).
\end{equation}
Then, estimates \eqref{estrho} and \eqref{acopmunu} hold in $\overline{\om}$.
\end{lemma}
  
{\begin{remark}\label{remlemafund}
Condition \eqref{reco1} means that $u_n$ is a recovery sequence in $\om$ for the functional 
\beq\label{functn}
w\in W^{1,p}(\om)^M\mapsto \int_\om F_n(x,Dw)\,dx,
\eeq
with the Dirichlet condition $w-u_n\in W^{1,p}_0(\omega)^M$. {  Since $w=u_n$ clearly satisfies $w-u_n\in W^{1,p}_0(\omega)^M$, this makes $u_n$ a recovery sequence without imposing any boundary condition.} In particular, condition \eqref{reco1}  is fulfilled if  for a fixed $f\in W^{-1,p}(\om)^M$, $u_n$ satisfies
    \begin{equation*}
       \int_\om F_n(x,Du_n)\,dx= \min\left\{\int_\om F_n\big(x,D(u_n+v)\big)\,dx-\langle f,v\rangle:\ v\in W^{1,p}_0(\om)^M \right\}.
    \end{equation*}
{  Assuming the differentiability of $F_n$ with respect to the second variable, it follows that} $u_n$ satisfies the variational equation
    \begin{equation*}
        \int_\om D_\xi F_n(x,Du_n):Dv\,dx -\langle f,v\rangle = 0,\quad\forall\, v\in W^{1,p}_0(\om)^M,
    \end{equation*}
   i.e. $u_n$ is a solution of 
        \begin{equation*}
        -\Dive\big(D_\xi F_n(x,Du)\big) = f\quad\mbox{in }\om,
    \end{equation*}
    where no boundary condition is imposed.\par
Assumption \eqref{reco1}  allows us to take into account very general boundary conditions. For example, if $u_n$ is a recovery sequence for \eqref{functn} with (non necessarily homogeneous) Dirichlet or Neumann boundary condition, then it also satisfies \eqref{reco1}.
 \end{remark}
\begin{remark}\label{remlemafund2}
  Condition \eqref{reco1} is equivalent to  the asymptotic minimizer property satisfied by $u_n$:
\begin{equation*}%\label{asymin}
\int_\om F_n(x,Du_n)\,dx \leq \int_\om F_n(x,Dw_n)\,dx+O_n,
\quad\forall\,w_n\;\;\hbox{with}\;\;w_n-u_n\rightharpoonup 0\ \hbox{ in } W^{1,p}_0(\om)^M.
\end{equation*}
We can check that if $u_n$ satisfies this condition in $\om$, then $u_n$ satisfies it in any open subset $\hat\om\subset\om$. To this end, it is enough to consider for a sequence $\hat w_n$ with $\hat w_n-u_n\in W^{1,p}_0(\hat\om)^M$, the extension
$$w_n:=\left\{\ba{ll} \hat w_n &\hbox{ in }\hat\om
\\
\dis u_n &\hbox{ in }\om\setminus\hat\om.\ea\right.$$
\end{remark}}
 
\begin{corollary}\label{clemafund}
  Let $F_n:\Omega\times \mathbb{R}^{M\times N}\to [0,\infty)$ be a sequence of Carath\'{e}odory functions satisfying \eqref{hip1} to \eqref{hipo3}. Consider two open sets $\omega_1,\omega_2\subset\Omega$ such that $\omega_1\cap\omega_2\neq {\rm\O}$, a sequence $u_n$ converging weakly in $W^{1,p}(\omega_1)^M$ to a function $u$ and a sequence $v_n$ converging weakly in $W^{1,p}(\omega_2)^M$ to a function $v$, such that
\begin{equation*}%}\label{reguuv}
u,v\in \left\{\ba{ll}\dis W^{1,\frac{pr}{r-1}}(\om_1\cap\om_2)^M, & \hbox{ if }1< p\leq N\!-\!1
\\ \ecart
\dis C^1(\om_1\cap\om_2)^M, & \hbox{ if }p>N\!-\!1,\ea\right.
\end{equation*}

  \begin{equation*}
    |Du_n|^p\stackrel{*}\rightharpoonup\varrho,\quad F_n({\cdot},Du_n)\stackrel{*}\rightharpoonup \mu\quad\mbox{ in }\M(\omega_1),
  \end{equation*}
  \begin{equation*}
    |Dv_n|^p\stackrel{*}\rightharpoonup\varpi,\quad F_n({\cdot},Dv_n)\stackrel{*}\rightharpoonup \nu\quad\mbox{ in }\M(\omega_2),
  \end{equation*}
  \begin{equation*}
    \exists \lim_{n\to\infty}\int_{\omega_1} F_n(x,Du_n)\,dx=\min\left\{\liminf_{n\to\infty} \int_{\omega_1} F_n(x,Dw_n)\,dx:\ w_n-u_n\rightharpoonup 0\mbox{ in }W^{1,p}_0(\omega_1)^M\right\},
  \end{equation*}
   \begin{equation*}
    \exists \lim_{n\to\infty}\int_{\omega_2} F_n(x,Dv_n)\,dx=\min\left\{\liminf_{n\to\infty} \int_{\omega_2} F_n(x,Dw_n)\,dx:\ w_n-v_n\rightharpoonup 0\mbox{ in }W^{1,p}_0(\omega_2)^M\right\}.
  \end{equation*}
Then, we have 
\begin{equation}\label{acomu-nu}
\ba{l}
|\mu-\nu|\leq
\\ \ecart
\begin{cases}
\dis C\big({h}^L\hskip-1pt+\hskip-1pt\mu^L\hskip-1pt+\hskip-1pt\nu^L\hskip-1pt+\hskip-1pt\varrho^L\hskip-1pt+\hskip-1pt\varpi^L\hskip-1pt+\hskip-1pt(1+({\textsc{a}}^L)^\frac{1}{r})|D(u-v)|^p\big)^\frac{p-1}{p}\hskip-1pt({\textsc{a}}^L)^\frac{1}{pr}\hskip-1pt|D(u-v)|\hbox{ a.e. in }\om_1\cap\om_2, & \hbox{if }1 \hskip-1pt <\hskip-1pt p\leq\hskip-1pt N \hskip-2pt -\hskip-2pt 1
\\ \ecart
\dis C\big(1+{h} ^{\textsc{a}}+\mu ^{\textsc{a}}+\nu ^{\textsc{a}}+\varrho ^{\textsc{a}}+\varpi ^{\textsc{a}}+|D(u-v)|^p\big)^\frac{p-1}{p}{\textsc{a}}\,|D(u-v)|\ {\textsc{a}}\hbox{-a.e. in }\om_1\cap\om_2, & \hbox{if }p>\hskip-1pt N \hskip-2pt -\hskip-2pt 1.
\end{cases}
\ea
\end{equation}
\end{corollary}

Lemma~\ref{lemafund} is itself based on the following compactness result.

\begin{lemma}\label{lem1}
Let $F_n:\Omega\times \mathbb{R}^{M\times N}\to [0,\infty)$ be a sequence of Carath\'{e}odory functions satisfying \eqref{hip1} to \eqref{hipo3}, and let $\om$ be an open subset of $\Omega$.
Consider a sequence $\{\xi_n\}\subset L^p(\om)^{M\times N}$ such that
\begin{equation}\label{codxinp}
F_n({\cdot},\xi_n)\stackrel{*}\rightharpoonup \Lambda\quad\mbox{and}\quad
|\xi_n|^p\stackrel{*}\rightharpoonup \Xi\quad\mbox{in }\M(\om).
\end{equation}
\begin{itemize}
    \item If $1<p\leq N\!-\!1$ and the sequence $\{\rho_n\}$ converges strongly to $\rho$ in $L^{\frac{pr}{r-1}}(\om)^{M\times N}$, then there exist a subsequence of $n$ and a function $\vartheta\in L^1(\om)$ such that
        \begin{equation}\label{cofudif}
          F_n({\cdot},\xi_n+\rho_n)-F_n({\cdot},\xi_n)\rightharpoonup \vartheta\quad\mbox{weakly in } L^1(\om),
        \end{equation}
        where $\vartheta$ satisfies
        \begin{equation}\label{acobvt}
          |\vartheta|\leq C\big({h}^L+\Lambda^L+\Xi^L+(1+({\textsc{a}}^L)^{\frac{1}{r}})|\rho|^p\big)^{\frac{p-1}{p}}({\textsc{a}}^L)^{\frac{1}{pr}}|\rho|\quad\mbox{a.e. in }\om.
        \end{equation}
    \item If $p>N\!-\!1$ and the sequence $\{\rho_n\}$ converges strongly to $\rho$ in $C^0(\overline{\om})^{M\times N}$, then there exist a subsequence of $n$ and a function $\vartheta \in L^1_{\textsc{a}}(\om)$ such that
        \begin{equation*}%\label{cofudifm}
          F_n({\cdot},\xi_n+\rho_n)\stackrel{*}\rightharpoonup\Lambda+\vartheta\,{\textsc{a}}\quad\mbox{in } \M(\om),
        \end{equation*}
        where $\vartheta$ satisfies
        \begin{equation}\label{acobvtb}
          |\vartheta|\leq C\big(1+{h}^{\textsc{a}}+\Lambda^{{\textsc{a}}}+\Xi^{\textsc{a}}+|\rho|^p\big)^{\frac{p-1}{p}}|\rho|
          \quad\mbox{${\textsc{a}}$-a.e. in }\om.
        \end{equation}
\end{itemize}
\end{lemma}

\subsection{Examples}\label{hyperelas}
In this section we give three examples of functionals $\F_n$ satisfying the assumptions \eqref{hip1} to \eqref{hipo3} of Theorem~\ref{thm}.
\begin{enumerate}
\item The first example illuminates the Lipschitz estimate \eqref{hipo2}. It is also based on a functional coercivity of type \eqref{hip1} rather than a pointwise coercivity.
\item The second example deals with the Saint Venant-Kirchhoff hyper-elastic energy (see, {\em e.g.}, \cite{Cia} Chapter~4).
\item The third example deals with an Ogden's type hyper-elastic energy (see, {\em e.g.}, \cite{Cia} Chapter~4).
\end{enumerate}
Let $\Om$ be a bounded set of $\RR^N$, $N\geq 2$. We denote for any function $u:\Om\to\RR^N$,
\beq\label{eEC}
e(u):={1\over 2}\left(Du+Du^T\right),\quad E(u):={1\over 2}\left(Du+Du^T+Du^TDu\right),
\quad C(u):=(I_N+Du)^T(I_N+Du).
\eeq

\subsubsection*{Example~1}
{  Let $p\in (1,\infty)$, and let $A_n$ be a symmetric tensor-valued function in $L^\infty\big(\Om;\L(\RR^{N\times N}_s)\big)$.}
We consider the energy density function defined by
\begin{equation*}%\label{exa1}
 {F_n(x,\xi):=\big|A_n(x)\xi^s:\xi^s\big|^\frac{p}{2}\quad\mbox{a.e. }x\in\Omega,\ \forall\,\xi\in\RR^{N\times N}.}
\end{equation*}
We assume that there exists $\al>0$ such that
\beq\label{cAn}
 {A_n(x)\xi:\xi\geq\al\,|\xi|^2,\quad\mbox{a.e. }x\in\Om,\ \forall\,\xi\in\RR^{N\times N}_s,}
\eeq
and that
{\beq\label{bAn}
\mbox{$|A_n|^{p\over 2}$ is bounded in $L^r(\Om)$ with }r\mbox{ defined by }\eqref{acoan}.
\eeq}
Then, the density $F_n$ and the associated functional
\begin{equation*}%\label{Fnexa1}
 {\F_n(u):=\int_\Om\big|A_ne(u):e(u)\big|^\frac{p}{2}dx}
\quad\mbox{for }u\in W^{1,p}_0(\Om)^N,
\end{equation*}
satisfy the conditions \eqref{hip1} to \eqref{hipo3} of Theorem~\ref{thm}.
\begin{proof}
Using successively \eqref{cAn} and the Korn inequality in $W^{1,p}_0(\Om)^N$ for $p>1$ (see, {\em e.g.}, \cite{Tem}), we have for any $u\in W^{1,p}_0(\Om)^N$,
\[
\F_n(u)=\int_\Om\big|A_ne(u):e(u)\big|^{\frac{p}{2}}\,dx\geq\al\int_\Om |e(u)|^p\,dx\geq\al\,C\int_\Om|Du|^p\,dx,
\]
which implies \eqref{hip1}. Conditions \eqref{acoL1} and \eqref{hipo3} are immediate.
It remains to prove condition~\eqref{hipo2} with estimate~\eqref{acoan}. Taking into account
that
\[
|D_\xi F_n(x,\xi)|=p\,\big|(A_n(x)\xi^s:\xi^s)^{p-2\over 2}A_n(x)\xi^s\big|\leq p\,\big|A_n(x)\xi^s:\xi^s|^{p-1\over 2}|A_n(x)|^{1\over 2},\quad\forall\,\xi\in \RR^{N\times N},\ \hbox{a.e. }x\in\Om,
\]
then using the mean value theorem and H\"older's inequality, we get

\[
\ba{ll}
\dis \big|F_n(x,\xi)-F_n(x,\eta)\big| & \dis  \le p\left((A_n\xi^s:\xi^s)^{1\over 2}+(A_n\eta^s:\eta^s)^{1\over 2}\right)^{p-1}|A_n|^{1\over 2}|\xi^s-\eta^s|
\\ \ecart
& \dis \leq p\,2^{(p-1)^2\over p}\big(F_n(x,\xi)+F_n(x,\eta)\big)^{p-1\over p}|A_n|^{1\over 2}|\xi-\eta|,
\ea
\]
for every $\xi,\eta\in\RR^{N\times N}$ and a.e. $x\in\Om$. This implies estimate \eqref{hipo2} with $h_n=0$ and $a_n=|A_n|^{p\over 2}$ bounded in $L^r(\Om)$.
\end{proof}
The two next examples belong to the class of hyper-elastic materials (see, {\em e.g.}, \cite{Cia}, Chapter~4).
\subsubsection*{Example~2}
For $N=3$, we consider the Saint Venant-Kirchhoff energy density defined by
\begin{equation}\label{exa2}
F_n(x,\xi):={\la_n(x)\over 2}\,\big[{\rm tr}\big(\tilde E(\xi)\big)\big]^2+\mu_n(x)\,\big|\tilde E(\xi)\big|^2,
\quad\mbox{a.e. }x\in\Omega,\ \forall\,\xi\in\RR^{3\times 3},
\end{equation}
where $\tilde E(\xi):={1\over 2}\left(\xi+\xi^T+\xi^T\xi\right)$, and $\la_n,\mu_n$ are the Lam\'e coefficients.
\par
We assume that there exists a constant $C>1$ such that
\beq\label{lamun}
\la_n,\mu_n\geq 0\;\;\mbox{a.e. in }\Om,\quad \essinf_\Om \left(\la_n+\mu_n\right)>C^{-1},\quad
\int_\Om(\la_n+\mu_n)\,dx\leq C.
\eeq
Then, the density $F_n$ and the associated functional (see definition \eqref{eEC})
\beq\label{Fnexa2}
\F_n(u):=\int_\Om\left({\la_n\over 2}\,\big[{\rm tr}\big(E(u)\big)\big]^2+\mu_n\,\big|E(u)\big|^2\right)dx
\quad\mbox{for }u\in W^{1,4}_0(\Om)^3,
\eeq
satisfy the conditions \eqref{hip1} to \eqref{hipo3} of Theorem~\ref{thm}.
\begin{proof}
There exists a constant $C>1$ such that we have for a.e. $x\in\Om$ and any $\xi\in\RR^{3\times 3}$,
\beq\label{cExi}
C^{-1}(\la_n+\mu_n)\,|\xi|^4-C\,(\la_n+\mu_n)\leq F_n(x,\xi)\leq C\,(\la_n+\mu_n)\,|\xi|^4+C\,(\la_n+\mu_n).
\eeq
Hence, we deduce that for a.e. $x\in\Om$ and any $\xi,\eta\in\RR^{3\times 3}$,
\[
\ba{ll}
\dis \big|F_n(x,\xi)-F_n(x,\eta)\big| & \dis \leq C\,(\la_n+\mu_n)\left(1+|\xi|^2+|\eta|^2\right)^{3\over 2}|\xi-\eta|
\\ \ecart
& \dis = C\left((\la_n+\mu_n)^{1\over 2}+(\la_n+\mu_n)^{1\over 2}|\xi|^2+(\la_n+\mu_n)^{1\over 2}|\eta|^2\right)^{3\over 2}(\la_n+\mu_n)^{1\over 4}|\xi-\eta|
\\ \ecart
& \dis \leq C\left((\la_n+\mu_n)^{1\over 2}+F_n(x,\xi)^{1\over 2}+F_n(x,\eta)^{1\over 2}\right)^{3\over 2}(\la_n+\mu_n)^{1\over 4}|\xi-\eta|
\\ \ecart
& \dis \leq C\,\big(\la_n+\mu_n+F_n(x,\xi)+F_n(x,\eta)\big)^{3\over 4}(\la_n+\mu_n)^{1\over 4}\,|\xi-\eta|,
\ea
\]
which implies estimate \eqref{hipo2} with $p=4$ and $h_n=a_n=\la_n+\mu_n$, while \eqref{acohn} and \eqref{acoan} are a straightforward consequence of \eqref{lamun}.
Moreover, by the first inequality of \eqref{cExi} combined with \eqref{lamun} we get that the functional~\eqref{Fnexa2} satisfies the coercivity condition \eqref{hip1}. Condition \eqref{acoL1} is immediate.
Finally, since we have
\[
\big[{\rm tr}\big(\tilde E(\la\xi)\big)\big]^2+\big|\tilde E(\la\xi)\big|^2\leq C\big(1+|\xi|^4\big),
\quad\forall\,\lambda\in [0,1],\ \forall\,\xi\in\mathbb{R}^{3\times 3},
\]
condition \eqref{hipo3} follows from the first inequality of \eqref{cExi}, which concludes the proof of the second example.
\end{proof}
\begin{remark}\label{rem.SVK}
The default of the Saint Venant-Kirchhoff model is that the function $F_n(x,\cdot)$ of \eqref{exa2} is not {\em polyconvex} (see \cite{Rao}). Hence,  we do not know if it is {\em quasiconvex}, or equivalently, if the functional $\F_n$ of \eqref{Fnexa2} is lower semi-continuous for the weak topology of $W^{1,4}(\Om)^3$ (see, {\em e.g.} \cite{Dac}, Chapter~4, for the notions of polyconvexity and quasiconvexity).
\end{remark}
\subsubsection*{Example~3}
For $N=3$ and $p\in[2,\infty)$, we consider the Ogden's type energy density defined by
{ 
\begin{equation}\label{exa3}
F_n(x,\xi):=a_n(x)\Big[{\rm tr}\big(\tilde C(\xi)^{p\over 2}-I_3\big)\Big]^+\quad\mbox{a.e. }x\in\Omega,\ \forall\,\xi\in\RR^{3\times 3},
\end{equation}
where $\tilde C(\xi):=(I_3+\xi)^T(I_3+\xi)$, and $t^+:=\max\,(t,0)$ for $t\in\RR$.
}
We assume that there exists a constant $C>1$ such that
\beq\label{an}
\essinf_{\Om}a_n>C^{-1}\quad\mbox{and}\quad\int_\Om a_n^r\,dx\leq C\quad
\mbox{with }
    \left\{\ba{cl}
       r>1, & \mbox{if }p=2
        \\ \ecart
       r=1, & \mbox{if } p > 2.
    \ea\right.
\eeq
Then, the density $F_n$ and the associated functional (see definition \eqref{eEC})
\beq\label{Fnexa3}
{ 
\F_n(u):=\int_\Om a_n(x)\Big[{\rm tr}\big(C(u)^{p\over 2}-I_3\big)\Big]^+dx\quad\mbox{for }u\in W^{1,p}_0(\Om)^3,
}
\eeq
satisfy the conditions \eqref{hip1} to \eqref{hipo3} of Theorem~\ref{thm}.
\begin{proof}
There exists a constant $C>1$ such that we have for a.e. $x\in\Om$ and any $\xi\in\RR^{3\times3}$,
\beq\label{cCxi}
C^{-1}a_n\,|\xi|^p-C\,a_n\leq F_n(x,\xi)\leq C\,a_n\,|\xi|^p+C\,a_n.
\eeq
This combined with the fact that the (well-ordered) eigenvalues of a symmetric matrix are Lipschitz functions (see, {\em e.g.}, \cite{Cia2}, Theorem~2.3-2), {implies} that for a.e. $x\in\Om$ and any $\xi,\eta\in\RR^N$, {we have}
\[
\ba{ll}
\dis \big|F_n(x,\xi)-F_n(x,\eta)\big| & \dis \leq C\,a_n(1+|\xi|+|\eta|)^{p-1}|\xi-\eta|
\\ \ecart
& \dis \leq C\,\big(a_n+a_n|\xi|^p+a_n|\eta|^p\big)^{p-1\over p}a_n^{1\over p}\,|\xi-\eta|
\\ \ecart
& \dis \leq C\,\big(a_n+F_n(x,\xi)+F_n(x,\eta)\big)^{p-1\over p}a_n^{1\over p}\,|\xi-\eta|,
\ea
\]
which implies estimate \eqref{hipo2} with $h_n=a_n$, while \eqref{acohn} and \eqref{acoan} are a straightforward consequence of \eqref{an}.
Moreover, by the first inequality of \eqref{cCxi} combined with \eqref{an} we get that the functional \eqref{Fnexa3} satisfies the coercivity condition \eqref{hip1}. Condition \eqref{acoL1} is immediate.
Finally, since we have
\[
{\rm tr}\big(\tilde C(\la\xi)^{p\over 2}\big)\leq C\big(1+|\xi|^p\big),
\quad\forall\,\lambda\in [0,1],\ \forall\,\xi\in\mathbb{R}^{3\times 3},
\]
condition \eqref{hipo3} follows from the first inequality of \eqref{cCxi}, which concludes the proof of the third example.
\end{proof}
\begin{remark}\label{rem.Ogden}
Contrary to Example~2, the function $F_n(x,\cdot)$ of \eqref{exa3} is polyconvex since it is the composition of the Ogden density energy defined for a.e. $x\in\Om$, by
\beq\label{Ogden}
{  
W_n(x,\xi):=a_n(x)\Big[{\rm tr}\big(\tilde C(\xi)^{p\over 2}-I_3\big)\Big]^+\quad\mbox{for }\xi\in\RR^{3\times 3},
}
\eeq
which is known to be polyconvex (see \cite{Bal}), by the non-decreasing convex function {$t\mapsto\,t^+$}.
However, in contrast with \eqref{Ogden} the function \eqref{exa3} does attain its minimum at $\xi=0$, namely in the absence of strain.
\end{remark}

\section{Proof of the results}
\subsection{Proof of the main results}
\noindent
{\bf Proof of Theorem~\ref{thm}.}
The proof is divided into two steps. In the first step we construct the limit functional $F$ and we prove the properties \eqref{proF0}, \eqref{proF1}, \eqref{proF2} satisfied by the function $F$. The second step is devoted to convergence \eqref{concthp}.
\par\smallskip\noindent
{\it First step:} Construction of $F$.
\par\noindent
Let $\F_n:W^{1,p}(\Om)^M\to[0,\infty]$ be the functional defined by
\begin{equation*}%\label{Fn}
\F_n(v)=\int_\Omega F_n(x,Dv)\,dx\quad\mbox{for }v\in W^{1,p}(\Omega)^M.
\end{equation*}
By the compactness $\Gamma$-convergence theorem (see {\em e.g.} \cite{Dal}, Theorem~8.5), there exists a subsequence of $n$, still denoted by $n$, such that $\F_n$ $\Gamma$-converges for the strong topology of $L^p(\Om)^M$ to a functional $\F:W^{1,p}(\Omega)^M\to~[0,\infty]$ with domain $\D(\F)$.
\par
Let $\xi$ be a matrix of a countable dense subset $D$ of $\RR^{M\times N}$ with $0\in D$. Since the linear function $x\mapsto\xi x$ belongs to $\D(\F)$ by \eqref{acotFn}, up to the extraction of a new subsequence, for any $\xi\in D$ there exists a recovery sequence $w_n^\xi$ in $W^{1,p}(\Omega)^M$ which converges strongly to $\xi x$ in $L^p(\Om)^M$ and such that
\begin{equation*}%\label{conFnwnxi}
F_n(\cdot,Dw_n^\xi)\stackrel{*}{\rightharpoonup}\mu^\xi\quad\mbox{and}\quad |Dw_n^\xi|^p\stackrel{*}{\rightharpoonup}\varrho^\xi
\quad\mbox{in }\M(\Om).
\end{equation*}
In particular, since $F_n(\cdot,0)=0$ we have $\mu^0=0$. Moreover, by estimates \eqref{estrho} and \eqref{acomu-nu} we have for any $\xi,\eta\in D$,

\beq\label{estrhoxi}
\varrho^\xi\leq\left\{\ba{ll}
C\big(|\xi|^p+|\xi|^p({\textsc{a}}^L)^{1\over r}+{h}+\mu^\xi+{\textsc{a}}^L\big)\quad\mbox{a.e. in }\omega, & \mbox{if }1<p\leq N\!-\!1
\\ \ecart
C\big(|\xi|^p {\textsc{a}}+{h}+\mu^\xi+{\textsc{a}}\big)\quad \textsc{a}\mbox{-a.e. in }\omega, & \mbox{if }p>N\!-\!1,
\ea\right.
\eeq

\begin{equation}\label{acomu-nuxi}
\ba{l}  |\mu^\xi-\mu^\eta|\leq\\ \ecart\dis
\left\{\hskip-6pt\ba{ll}\dis C\big({h}^L\hskip-1pt+\hskip-1pt(\mu^\xi)^L\hskip-1pt+\hskip-1pt(\mu^\eta)^L\hskip-1pt+\hskip-1pt(\varrho^\xi)^L\hskip-1pt+\hskip-1pt(\rho^\eta)^L\hskip-1pt+\hskip-1pt(1+({\textsc{a}}^L)^\frac{1}{r})|\xi-\eta|^p\big)^\frac{p-1}{p}\hskip-1pt({\textsc{a}}^L)^\frac{1}{pr}\hskip-1pt|\xi-\eta|\hbox{ a.e. in }\Om, & \hbox{if }1<p\leq\hskip-1pt N \hskip-2pt -\hskip-2pt 1
\\ \ecart
\dis C\left(1+{h} ^{\textsc{a}}+(\mu^\xi) ^{\textsc{a}}+(\mu^\eta) ^{\textsc{a}}+(\varrho^\xi) ^{\textsc{a}}+(\varrho^\eta) ^{\textsc{a}}+|\xi-\eta)|^p\right)^\frac{p-1}{p}{\textsc{a}}\,|\xi-\eta|\ {\textsc{a}}\hbox{-a.e. in }\Om, & \hbox{if }p>\hskip-1pt N \hskip-2pt -\hskip-2pt 1.\ea\right.\ea\eeq

Hence, by a continuity argument we can define a function $F:\Om\times\RR^{M\times N}\to [0,\infty)$ satisfying \eqref{proF0}, \eqref{proF2} and such that
\beq\label{F}
\mu^\xi=\left\{\ba{ll}
F(\cdot,\xi), & \hbox{if }1<p\leq N\!-\!1
\\ \ecart
F(\cdot,\xi)\,{\textsc{a}}, & \mbox{if }p>N\!-\!1,
\ea\right.
\quad\forall\,\xi\in D,
\eeq
where the property \eqref{proF1} is deduced from \eqref{estrhoxi}, \eqref{acomu-nuxi}.
\par\medskip\noindent
{\it Second step:} Proof of convergence \eqref{concthp}.
\par\noindent
Let $\om$ be an open set of $\Om$, let $\{u_n\}$ be a sequence fulfilling \eqref{reco1bis}, which converges weakly in $W^{1,p}(\om)^M$ to a function $u$ satisfying \eqref{regulu}, and let $\xi\in D$. Since $F_n(\cdot,Du_n)$ is bounded in $L^1(\Om)$, there exists a subsequence of $n$, still denoted by $n$, such that
\beq\label{conFnDun}
F_n(\cdot,Du_n)\stackrel{*}{\rightharpoonup}\mu\quad\mbox{and}\quad |Du_n|^p\stackrel{*}{\rightharpoonup}\varrho\quad\mbox{in }\M(\Om).
\eeq
Applying Corollary~\ref{clemafund} to the sequences $u_n$ and $v_n=w_n^\xi$, we have
\begin{equation*}
\ba{l}  |\mu-\mu^\xi|\leq\\ \ecart\dis
\left\{\hskip-6pt\ba{ll}\dis C\big({h}^L\hskip-1pt+\hskip-1pt\mu^L\hskip-1pt+\hskip-1pt(\mu^\xi)^L\hskip-1pt+\hskip-1pt\varrho^L\hskip-1pt+\hskip-1pt(\varrho^\xi)^L\hskip-1pt+\hskip-1pt(1+({\textsc{a}}^L)^\frac{1}{r})|Du-\xi|^p\big)^\frac{p-1}{p}\hskip-1pt({\textsc{a}}^L)^\frac{1}{pr}\hskip-1pt|Du-\xi|\hbox{ a.e. in }\om, &\hbox{if }1<p\leq\hskip-1pt N \hskip-2pt -\hskip-2pt 1
\\ \ecart\dis C\big(1+{h} ^{\textsc{a}}+\mu ^{\textsc{a}}+(\mu^\xi) ^{\textsc{a}}+\varrho ^{\textsc{a}}+(\varrho^\xi) ^{\textsc{a}}+|Du-\xi|^p\big)^\frac{p-1}{p}{\textsc{a}}\,|Du-\xi|\quad{\textsc{a}}\hbox{-a.e. in }\om, &\hbox{if }p>\hskip-1pt N \hskip-2pt -\hskip-2pt 1.\ea\right.\ea
\end{equation*}

Using \eqref{estrhoxi}, \eqref{F} and the continuity of $F(x,\xi)$ with respect to $\xi$, we get that
\beq\label{muFDu}
\mu=\left\{\ba{ll}
F(\cdot,Du), & \hbox{if }1<p\leq N\!-\!1
\\ \ecart
F(\cdot,Du)\,{\textsc{a}}, & \mbox{if }p>N\!-\!1.
\ea\right.
\eeq
Note that since the limit $\mu$ is completely determined by $F$, the first convergence of \eqref{conFnDun} holds for the whole sequence,
which concludes the proof.
\cqfd
\par\bigskip\noindent
{\bf Proof of Theorem~\ref{thm2}.}
The proof is divided into two steps.
\par\ms\noindent
{\it First step:} The case where $V=\{\hat{u}\}+W^{1,p}_0(\om)^M$.
\par\noindent
Fix a function $\hat{u}$ satisfying \eqref{regulu1}, and define the set $V:=\{\hat{u}\}+W^{1,p}_0(\om)^M$.
Let $u\in V$ such that
$$
u\in \left\{\ba{ll}\dis W^{1,\frac{pr}{r-1}}(\om)^M, & \hbox{ if }1< p\leq N\!-\!1
\\ \ecart
\dis C^1(\overline{\om})^M, &\hbox{ if }p>N\!-\!1.
\ea\right.
$$
which is extended by $\hat{u}$ in $\Om\setminus\om$, and consider a recovery sequence $\{u_n\}$ for $\F_n^V$ of limit $u$.
There exists a subsequence of $n$, still denoted by $n$, such that the first convergences of \eqref{comedu2} and \eqref{devr2} hold.
By Theorem~\ref{thm} convergences \eqref{concthp} are satisfied in $\om$, which implies \eqref{muFDu}.
Now, applying the estimate \eqref{acopmunu} of Lemma~\ref{lemafund2} with $u_n$ and $v_n=u$, it follows that
\[
\mu\leq\nu\;\;\mbox{in }\overline{\om}\quad\mbox{with}\quad F_n(\cdot,Dv_n)\stackrel{*}\rightharpoonup\nu\;\;\mbox{in }\M(\Om),
\]
where the convergence holds up to a subsequence.
Then, using estimate \eqref{hipo2} with $\eta=0$ and H\"older's inequality, we have for any $\ph\in L^\infty\big(\Om;[0,1]\big)$ with compact support in $\Om$,
\[
\ba{l}
\dis \int_\Om\ph\,F_n(x,Du)\,dx\leq
\\ \ecart\left\{\ba{ll}
\dis \left(\into\ph\,\big(h_n+F_n(x,Du)+|Du|^p\big)\,dx\right)^{p-1\over p}\left(\into\ph\,a_n^r\,dx\right)^{1\over pr}\left(\ph\,|Du|^{pr\over r-1}dx\right)^{r-1\over pr}, & \mbox{if }1<p\leq N\!-\!1
\\ \ecart
\dis \left(\into\ph\,\big(h_n+F_n(x,Du)+|Du|^p\big)\,dx\right)^{p-1\over p}\left(\into\ph\,a_n\,dx\right)^{1\over p}\|Du\|_{L^\infty(\Om)^M},
& \mbox{if }p>N\!-\!1,
\ea\right.\ea
\]
which implies that $\nu$ is absolutely continuous with respect to the Lebesgue measure if $1<p\leq N\!-\!1$, and absolutely continuous with respect to measure ${\textsc{a}}$ if $p>N\!-\!1$. Due to condition \eqref{omega} in both cases the equality $\nu(\partial\om)=0$ holds, so does with $\mu$.
This combined with \eqref{concthp} and \eqref{muFDu} yields
\begin{equation*}%\label{muFDu2}
\lim_{n\to\infty}\int_\om F_n(x,Du_n)\,dx=
\left\{\ba{ll}
\dis \int_\om F(x,Du)\,dx, & \hbox{if }1<p\leq N\!-\!1
\\ \ecart
\dis \int_\om F(x,Du)\,d{\textsc{a}}, & \mbox{if }p>N\!-\!1,
\ea\right.
\end{equation*}
which concludes the first step.
\par\ms\noindent
{\it Second step:} The general case.
\par\noindent
Let $V$ be a subset of $W^{1,p}(\om)^M$ satisfying \eqref{V}.
Let $u$ be a function such that
\[
u\in
\left\{\ba{lll}
V\cap W^{1,{pr\over r-1}}(\Om)^M, & \mbox{if }1<p\leq N\!-\!1
\\ \ecart
V\cap C^1(\overline{\Om})^M, & \mbox{if }p>N\!-\!1,
\ea\right.
\]
 and define the set $\tilde{V}:=\{u\}+W^{1,p}_0(\om)^M$.
Consider a recovery sequence $\{u_n\}$ for $\F_n^V$ given by \eqref{FnV} of limit $u$, and a recovery sequence $\{\tilde{u}_n\}$ for $\F_n^{\tilde{V}}$ of limit $u$. By virtue of Theorem~\ref{thm} the convergences \eqref{concthp} hold for both sequences $\{u_n\}$ and $\{\tilde{u}_n\}$.
Hence, since $\om$ is an open set, {and $F_n(x,Du_n)$ is non-negative}, we have
\beq\label{concthp2}
\left.\ba{ll}
\hbox{if }1<p\leq N\!-\!1, & \dis \int_\om F(x,Du)\,dx
\\ \ecart
\hbox{if } p>N\!-\!1, & \dis \int_\om F(x,Du)\,d{\textsc{a}}
\ea\right\}
\leq \liminf_{n\to\infty}\int_\om F_n(x,Du_n)\,dx.
\eeq
Moreover, since $\tilde{u}_n-u_n\rightharpoonup 0$ in $W^{1,p}_0(\om)^M$, $\tilde{u}_n\in V$ by property \eqref{V} and because $\{u_n\}$ is a recovery sequence for $\F_n^V$, $\{\tilde{u}_n\}$ is an admissible sequence for the minimization problem \eqref{reco1bis}, which implies that
\beq\label{concthp3}
\exists\,\lim_{n\to\infty}\int_\om F_n(x,Du_n)\,dx
\leq\liminf_{n\to\infty}\int_\om F_n(x,D\tilde{u}_n)\,dx.
\eeq
On the other hand, by the first step applied with $\tilde{u}=u$ and the set $\tilde{V}$, we have
\beq\label{concthp4}
\lim_{n\to\infty}\int_\om F_n(x,D\tilde{u}_n)\,dx=
\left\{\ba{ll}
\dis \int_\om F(x,Du)\,dx, & \hbox{if }1<p\leq N\!-\!1
\\ \ecart
\dis \int_\om F(x,Du)\,d{\textsc{a}}, & \hbox{if } p>N\!-\!1.
\ea\right.
\eeq
Therefore, combining \eqref{concthp2}, \eqref{concthp3}, \eqref{concthp4}, for the sequence $n$ obtained in Theorem~\ref{thm}, the sequence $\big\{\F_n^V\big\}$ $\Gamma$-converges to some functional $\F^V$ satisfying \eqref{FV} with $v=u$, which concludes  the proof of Theorem~\ref{thm2}.
\cqfd

\subsection{Proof of the lemmas}
\noindent
{\bf Proof of Lemma~\ref{lem1}.}
  Assume that $1<p\leq N\!-\!1$. Using \eqref{acotFnb}, we have
  \begin{equation*}
    F_n(x,\xi_n+\rho_n)\leq (p-1)h_n+(2p-1)F_n(x,\xi_n)+(p-1)\big(|\xi_n+\rho_n|^p+|\xi_n|^p\big)+a_n|\rho_n|^p\quad\mbox{a.e. in }\om.
  \end{equation*}
  From this we deduce that $\{F_n({\cdot},\xi_n+\rho_n)\}$ is bounded in $L^1(\om)$. Moreover, by \eqref{hipo2}, we have
  \begin{equation*}
    \big|F_n(x,\xi_n+\rho_n)-F_n(x,\xi_n)\big| \leq \big(h_n+F_n(x,\xi_n+\rho_n)+F_n(x,\xi_n)+|\xi_n+\rho_n|^p+|\xi_n|^p\big)^\frac{p-1}{p}a_n^\frac{1}{p}|\rho_n|\quad\mbox{a.e. in }\om,
  \end{equation*}
  where, thanks to the strong convergence of $\{\rho_n\}$ in $L^{\frac{pr}{r-1}}(\om)^{M\times N}$, we can show that the right-hand side is bounded in $L^1(\om)$ and equi-integrable. Indeed, taking into account 
  \[\frac{p-1}{p}+\frac{1}{pr}+\frac{r-1}{pr}=1, \]
  we have the boundedness in $L^1(\om)$, while the strong convergence of  $\{\rho_n\}$ in $L^{\frac{pr}{r-1}}(\om)^{M\times N}$ implies that $\{|\rho_n|^{\frac{pr}{r-1}}\}$ is equi-integrable and therefore, the equi-integrability of the right-hand side. By the Dunford-Pettis theorem, extracting a subsequence if necessary, we conclude \eqref{cofudif}, which, together with \eqref{codxinp}, in particular implies
\[F_n({\cdot},\xi_n+\rho_n)\stackrel{*}\rightharpoonup \Lambda + \vartheta\quad\mbox{in }\M(\om).\]

  Moreover, for any ball $B\subset \om$, we have
  \begin{equation*}
  \begin{aligned}
        &\int_B\big|F_n(x,\xi_n+\rho_n)-F_n(x,\xi_n)\big|\,dx
        \\
    \leq& \int_B \big(h_n+F_n(x,\xi_n+\rho_n)+F_n(x,\xi_n)+|\xi_n+\rho_n|^p+|\xi_n|^p\big)^\frac{p-1}{p}a_n^\frac{1}{p}|\rho_n|dx\\
    \leq& \left(\int_B \big(h_n+F_n(x,\xi_n+\rho_n)+F_n(x,\xi_n)+C|\xi_n|^p+C|\rho_n|^p\big)\,dx\right)^\frac{p-1}{p}\left(\int_B a_n^rdx\right)^\frac{1}{pr}\left(\int_B|\rho_n|^\frac{pr}{r-1}dx\right)^\frac{r-1}{pr},
\end{aligned}
\end{equation*}

  \noindent
  which, passing to the limit, implies
  \[\int_B |\vartheta|dx\leq \left(({h}+2\Lambda+\vartheta+C\Xi)(\overline{B})+C\int_B|\rho|^pdx\right)^{\frac{p-1}{p}}{\textsc{a}}(\overline{B})^{\frac{1}{pr}}\left(\int_B|\rho|^{\frac{pr}{r-1}}dx\right)^{\frac{r-1}{pr}},\]

  \noindent
  and then, dividing by $|B|$, the measures differentiation theorem shows that
  \begin{equation}\label{acovtt}
    |\vartheta|\leq \big({h}^L+2\Lambda^L+\vartheta+C\Xi+C|\rho|^p\big)^{\frac{p-1}{p}}({\textsc{a}}^L)^\frac{1}{pr}|\rho|\quad\mbox{ a.e. in }\om.
  \end{equation}

  Using Young's inequality in \eqref{acovtt}
  \[|\vartheta|\leq \frac{p-1}{p} \big({h}^L+2\Lambda^L+\vartheta+C\Xi^L+C|\rho|^p\big)+\frac{1}{p}({\textsc{a}}^L)^\frac{1}{r}|\rho|^p\quad\mbox{a.e. in }\om,\]

  \noindent
  and then
  \[|\vartheta|\leq C \big({h}^L+\Lambda^L+\Xi^L+(1+({\textsc{a}}^L)^\frac{1}{r})|\rho|^p\big)\quad\mbox{a.e. in }\om,\]

  \noindent
  which substituted in \eqref{acovtt} shows \eqref{acobvt}.\newline

  Assume now that $p>N\!-\!1$. Again, using \eqref{acotFnb} we deduce that $\{F_n({\cdot},\xi_n+\rho_n)\}$ is bounded in $L^1(\om)$, and thanks to \eqref{hipo2} we get
   \begin{equation*}%\label{acoFn+Fn}
   \big|F_n(x,\xi_n+\rho_n)-F_n(x,\xi_n)\big| \leq \big(h_n+F_n(x,\xi_n+\rho_n)+F_n(x,\xi_n)+|\xi_n+\rho_n|^p+|\xi_n|^p\big)^\frac{p-1}{p}a_n^\frac{1}{p}|\rho_n|\quad\mbox{a.e. in } \om.
  \end{equation*}

  \noindent
  Consequently, the sequence $\{F_n({\cdot},\xi_n+\rho_n)-F_n({\cdot},\xi_n)\}$ is bounded in $L^1(\om)$. Extracting a subsequence if necessary, the sequence $\{F_n({\cdot},\xi_n+\rho_n)-F_n({\cdot},\xi_n)\}$ weakly-$*$ converges in $\M(\om)$ to a measure $\Theta$, which, together with \eqref{codxinp}, implies
  \begin{equation*}
    F_n({\cdot},\xi_n+\rho_n)\stackrel{*}\rightharpoonup \Lambda+\Theta\quad\mbox{in }\M(\om).
  \end{equation*} Furthermore, if $E$ is a measurable subset of $\om$, then, using H\"older's inequality, we have
  \begin{equation*}
  \begin{aligned}
        & \int_E \big|F_n(x,\xi_n+\rho_n)-F_n(x,\xi_n)\big|\,dx
        \\
    \leq & \int_E \big(h_n+F_n(x,\xi_n+\rho_n)+F_n(x,\xi_n)+|\xi_n+\rho_n|^p+|\xi_n|^p\big)^\frac{p-1}{p}a_n^\frac{1}{p}|\rho_n|\,dx
    \\
    \leq & \|\rho_n\|_{L^\infty(\om)^{M\times N}}\left(C\|\rho_n\|_{L^\infty(\om)^{M\times N}}^p+\int_E \big(h_n+F_n(x,\xi_n+\rho_n)+F_n(x,\xi_n)+C|\xi_n|^p\big)\,dx\right)^\frac{p-1}{p}\left(\int_E a_ndx\right)^\frac{1}{p},
  \end{aligned}
  \end{equation*}

\noindent
which, passing to the limit, shows that $\Theta$ is absolutely continuous with respect to ${\textsc{a}}$. By the Radon-Nikodym theorem, there exists $\vartheta\in L^1_{\textsc{a}}({\om})$ such that
\[\Theta=\vartheta{\textsc{a}}\quad\mbox{in }\M(\om).\]
From the previous expression and using the measures differentiation theorem, we get \eqref{acobvtb}.
\cqfd
\par\bigskip\noindent
{\bf Proof of Lemma \ref{lemafund}}. Let $x_0\in \omega$ and two numbers $0<R_1<R_2$ with $B(x_0,R_2)\subset\omega$. Lemma 2.6 in \cite{BrCa5} gives the existence of a sequence of closed sets
\[
U_n\subset[R_1,R_2], \mbox{ with }|U_n|\geq \frac{1}{2}(R_2-R_1),
\]
such that defining
\begin{equation*}
\bar{u}_n(r,z)=u_n(x_0+rz),\quad \bar{u}(r,z)=u(x_0+rz),\quad r\in(0,R_2),\ z\in S_{N\!-\!1},
\end{equation*}
  we have
  \begin{equation}\label{conun}
  \|\bar{u}_n-\bar{u}\|_{C^0(U_n;X)}\to0,
  \end{equation}
where $X$ is the space defined by
 \begin{equation*}%\label{deX}
  X:=\begin{cases}
    L^s(S_{N\!-\!1})^M,\mbox{ with }1\leq s<\frac{(N\!-\!1)p}{N\!-\!1-p}, & \mbox{if }1<p<N\!-\!1,
    \\
    L^s(S_{N\!-\!1})^M,\mbox{ with }1\leq s<\infty, & \mbox{if }p=N\!-\!1,
    \\
    C^0(S_{N\!-\!1})^M, & \mbox{if }p>N\!-\!1.
    \end{cases}
  \end{equation*}
 \par\noindent
 For the rest of the prove we assume $1<p\leq N\!-\!1$ because the case $p>N\!-\!1$ is quite similar. \\
 
We define $\bar{\varphi}_n\in W^{1,\infty}(0,\infty)$ by
  \begin{equation}\label{dephir}
    \bar{\varphi}_n(r)=\begin{cases}
      1,&\mbox{if }0<r<R_1,
      \\
    \dis  \frac{1}{|U_n|}\int_r^{R_2}\chi_{U_n}ds, &\mbox{if }R_1<r<R_2,
      \\
      0,&\mbox{if }R_2<r,
    \end{cases}
  \end{equation}
and
  \begin{equation*}%\label{barphn}
  \varphi_n(x)=\bar{\varphi}_n(|x-x_0|).
  \end{equation*}  
 Applying the coercivity inequality \eqref{hip1} to the sequence $\ph_n(u_n-u)$ and using $F_n(\cdot,0)=0$, $\ph_n=1$ in $B(x_0,R_1)$,  we get
\[
{ 
 \ba{l}
 \dis  \al\int_{B(x_0,R_1)}|Du_n-Du|^p\,dx \leq \al\int_{B(x_0,R_2)}\big|D\big(\ph_n(u_n-u)\big)\big)\big|^p\,dx \\ \ecart
\dis\leq \int_{B(x_0,R_2)}F_n\big(x,D(\varphi_n(u_n-u))\big)\,dx = \int_{B(x_0,R_2)}F_n\big(x,\ph_nDu_n-\ph_nDu+(u_n-u)\otimes\nabla\ph_n\big)\,dx.
 \ea
 }
 \]
 By the convergence \eqref{cofudif} with $\xi_n:=\ph_nDu_n$, $\rho_n:=-\ph_nDu+(u_n-u)\otimes\nabla\ph_n$, and by estimate \eqref{hipo3} we obtain up to a subsequence
 \[
 \ba{l}
 \dis \lim_{n\to\infty}\int_{B(x_0,R_2)}F_n\big(x,\ph_nDu_n-\ph_nDu+(u_n-u)\otimes\nabla\ph_n\big)\,dx
 \\ \ecart
 \dis \leq\lim_{n\to\infty}\int_{B(x_0,R_2)}F_n\big(x,\ph_nDu_n)\,dx+\int_{B(x_0,R_2)}\vartheta\,dx
\\ \ecart
 \dis \leq C({h}+\mu)\big(\overline{B}(x_0,R_2)\big)+\int_{B(x_0,R_2)}\vartheta\,dx,
 \ea
 \]
 with
 \[
 |\vartheta|\leq C\big({h}^L+\mu^L+\varrho^L+(1+({\textsc{a}}^L)^{\frac{1}{r}})|Du|^p\big)^{\frac{p-1}{p}}({\textsc{a}}^L)^{\frac{1}{pr}}
 |Du|\mbox{ a.e. in }\om.
 \]
 Indeed, thanks to \eqref{conun} the sequence $(u_n-u)\otimes\nabla\varphi_n$ converges strongly to $0$ in $L^{pr\over r-1}(\om)^{M\times N}$ taking into account the inequality
\[
{(N-1) p\over N-1-p}\geq{pr\over r-1}.
\]
Hence, we deduce from the previous estimates that
\[
\ba{ll}
\varrho\big(B(x_0,R_1)\big)\leq &\dis C({h}+\mu)\big(\overline{B}(x_0,R_2)\big)+C\int_{B(x_0,R_1)}|D u|^p\,dx
\\ \ecart
&\dis +\,C\int_{B(x_0,R_2)}\left(\big({h}^L+\mu^L+\varrho^L+(1+({\textsc{a}}^L)^{\frac{1}{r}})|Du|^p\big)^{\frac{p-1}{p}}({\textsc{a}}^L)^{\frac{1}{pr}} |Du|\right)dx.
 \ea
\]
Taking $R_2$ such that
 \[
 ({h}+\mu)\big(\{|x-x_0|=R_2\}\big)=0,
 \]
 which holds true except for a countable set $E_{x_0}\subset \big(0, \dist(x_0,\partial\omega)\big)$, and making $R_1$ tend to $R_2$,  we get that
 \[
\ba{ll}
\varrho\big(B(x_0,R_2)\big)\leq & \dis C({h}+\mu)\big(B(x_0,R_2)\big)+C\int_{B(x_0,R_2)}|D u|^p\,dx
\\ \ecart
& \dis +\,C\int_{B(x_0,R_2)}\left(\big({h}^L+\mu^L+\varrho^L+(1+({\textsc{a}}^L)^{\frac{1}{r}})|Du|^p\big)^{\frac{p-1}{p}}({\textsc{a}}^L)^{\frac{1}{pr}} |Du|\right)dx,
 \ea
\]
for any $R_2\in\big(0, \dist(x_0,\partial\omega)\big)\setminus E_{x_0}$.
Then, by the measures differentiation theorem it follows that
\begin{equation*}%\label{estrho2}
\varrho\leq C\left(|D u|^p+{h}+\mu\right)+C\left(\big({h}^L+\mu^L+\varrho^L+(1+({\textsc{a}}^L)^{\frac{1}{r}})|Du|^p\big)^{\frac{p-1}{p}}\right)({\textsc{a}}^L)^{\frac{1}{pr}} |Du|.
\end{equation*}
Finally, the Young inequality yields the desired estimate \eqref{estrho}.

Now consider $\{u_n\}$ and $\{v_n\}$ as in the statement of the lemma. Let $x_0\in \omega$ and $0<R_0<R_1<R_2$ with $B(x_0,R_2)\subset\omega$. Again using Lemma 2.6 in \cite{BrCa5} there exist two sequences of closed sets
\[V_n\subset[R_0,R_1],\quad U_n\subset[R_1,R_2],\]

\noindent
  with
  \[|V_n|\geq \frac{1}{2}(R_1-R_0),\quad|U_n|\geq \frac{1}{2}(R_2-R_1),\]

  \noindent
  such that defining
  \begin{equation*}
  \begin{aligned}
    &\bar{u}_n(r,z)=u_n(x_0+rz),&&\quad\bar{v}_n(r,z)=v_n(x_0+rz),&\quad r\in(0,R_2),\ z\in S_{N\!-\!1},\\
    &\bar{u}(r,z)=u(x_0+rz),&&\quad\bar{v}(r,z)=v(x_0+rz),&\quad r\in(0,R_2),\ z\in S_{N\!-\!1},
  \end{aligned}
  \end{equation*}

  \noindent
  we have
  \begin{equation*}%\label{conunvn}
    \|\bar{u}_n-\bar{u}\|_{C^0(U_n;X)}\to0,\quad\|\bar{v}_n-\bar{v}\|_{C^0(V_n;X)}\to0.
  \end{equation*}
    Then, consider the function $\bar{\ph}_n$ defined by \eqref{dephir} and the function $\bar{\psi}_n\in W^{1,\infty}(0,\infty)$ defined by

  \begin{equation*}%\label{depsir}
    \bar{\psi}_n(r)=\begin{cases}
      1,&\mbox{if }0<r<R_0,\\ \dis
      \frac{1}{|V_n|}\int_r^{R_1}\chi_{V_n}ds, &\mbox{if }R_0<r<R_1,\\
      0,&\mbox{if }R_1<r.
    \end{cases}
  \end{equation*}

  \noindent
  From these sequences we define $w_n\in W^{1,p}(\omega)^M$ by
  \[w_n=\psi_n(v_n-v+u)+\varphi_n(1-\psi_n)u+(1-\varphi_n)u_n,\]

  \noindent
  with
  \begin{equation*}%\label{barphnpsin}
  \varphi_n(x)=\bar{\varphi}_n(|x-x_0|),\quad \psi_n(x)=\bar{\psi}_n(|x-x_0|),
  \end{equation*}  

  \noindent
  {\em i.e.}
  \beq\label{wn}
  w_n=\begin{cases}
    v_n-v+u,&\mbox{if } |x-x_0|<R_0,\\
    \psi_n(v_n-v)+u,&\mbox{if }R_0<|x-x_0|<R_1,\\
    \varphi_nu+(1-\varphi_n)u_n,&\mbox{if }R_1<|x-x_0|<R_2,\\
    u_n,&\mbox{if }R_2<|x-x_0|, x \in \omega.
  \end{cases}
  \eeq

It is clear that, for a subsequence, $w_n$ converges a.e. to $u$. Using then that $w_n-u_n$ is in $W^{1,p}_0(\om)^M$ and that, thanks to $\varphi_n$, $\psi_n$ bounded in $W^{1,\infty}(\Om)$, $w_n$ is bounded in $W^{1,p}(\om)^M$, we get 
  \begin{equation*}
    w_n-u_n \rightharpoonup 0\quad\mbox{weakly in }W^{1,p}_0(\omega).
  \end{equation*}

  Thus, from \eqref{reco1} we deduce
  \begin{equation*}
  \begin{aligned}
    \int_\omega F_n(x,Du_n)\,dx\leq&\int_\omega F_n(x,Dw_n)\,dx+O_n\\
    = &\int_{B(x_0,R_0)}F_n\big(x,D(v_n-v+u)\big)\,dx+\int_{\{R_2<|x-x_0|\}\cap \omega}F_n(x,Du_n)\,dx\\
    &+\int_{\{R_0<|x-x_0|<R_1\}}F_n\big(x,\psi_n D(v_n-v)+Du+(v_n-v)\otimes \nabla \psi_n\big)\,dx\\
    &+\int_{\{R_1<|x-x_0|<R_2\}} F_n\big(x,\varphi_nDu+(1-\varphi_n)Du_n+(u-u_n)\otimes \nabla\varphi_n\big)\,dx+O_n,
  \end{aligned}
  \end{equation*}
  what implies, in particular
  \begin{equation}\label{defu1}
  \begin{aligned}
    \int_{B(x_0,R_2)} F_n(x,Du_n)\,dx\leq&\int_{B(x_0,R_0)}F_n\big(x,D(v_n-v+u)\big)\,dx\\
    &+\int_{\{R_0<|x-x_0|<R_1\}}F_n\big(x,\psi_n D(v_n-v)+Du+(v_n-v)\otimes \nabla \psi_n\big)\,dx\\
    &+\int_{\{R_1<|x-x_0|<R_2\}} F_n\big(x,\varphi_nDu+(1-\varphi_n)Du_n+(u-u_n)\otimes \nabla\varphi_n\big)\,dx+O_n.
  \end{aligned}
  \end{equation}

  To estimate the first term on the right-hand side of this inequality, we use Lemma \ref{lem1} with $\xi_n=Dv_n$, $\rho_n=D(-v+u)$, which take into account \eqref{comedv},  gives
  \begin{equation}\label{esdef1}
  \begin{aligned}
    &\ \int_{B(x_0,R_0)} F_n\big(x,D(v_n-v+u)\big)\,dx \\
    \leq&\ \nu\big(\overline{B}(x_0,R_0)\big)+C\int_{B(x_0,R_0)} \big({h}^L+\nu^L+\varpi^L+(1+({\textsc{a}}^L)^{\frac{1}{r}})|D(u-v)|^p\big)^{\frac{p-1}{p}}({\textsc{a}}^L)^\frac{1}{pr}|D(u-v)|dx+O_n.
  \end{aligned}
  \end{equation}
   For the second term, we use again Lemma \ref{lem1} with $\xi_n=\psi_n Dv_n$ and $\rho_n = -\psi_n Dv+Du+(v_n-v)\otimes \nabla \psi_n$. Therefore, up to subsequence it holds  
  \begin{equation}\label{esdef2}
  \begin{aligned}
   &\ \int_{\{R_0<|x-x_0|<R_1\}}F_n\big(x,\psi_n D(v_n-v)+Du+(v_n-v)\otimes \nabla \psi_n\big)\,dx\\ \ecart
    \leq&\ C({h}+\nu+\varpi)\big(\{R_0\leq|x-x_0|\leq R_1\}\big)\\
    &+C\int_{\{R_0<|x-x_0|<R_1\}} \big({h}^L+\nu^L+\varpi^L+(1+({\textsc{a}}^L)^\frac{1}{r})(|Dv|^p+|Du|^p)\big)^\frac{p-1}{p}({\textsc{a}}^L)^\frac{1}{pr}(|Du|+|Dv|)\,dx+O_n.
  \end{aligned}
  \end{equation}
  The third term is analogously estimated by Lemma \ref{lem1} with $\xi_n = (1-\varphi_n) Du_n$ and $\rho_n = \varphi_n Du + (u-u_n)\otimes \nabla\varphi_n$. Extracting a subsequence if necessary, it yields
  \begin{equation}\label{esdef3}
  \begin{aligned}
     &\int_{\{R_1<|x-x_0|<R_2\}}F_n\big(x,\varphi_nDu+(1-\varphi_n)Du_n+(u-u_n)\otimes \nabla\varphi_n\big)\,dx\\ \ecart
     \leq&\ C({h}+\mu+\varrho)\big(\{R_1\leq|x-x_0|\leq R_2\}\big)\\	
    &+C\int_{\{R_1<|x-x_0|<R_2\}}\big({h}^L+\mu^L+\varrho^L+(1+({\textsc{a}}^L)^\frac{1}{r})|Du|^p\big)^\frac{p-1}{p}({\textsc{a}}^L)^\frac{1}{pr}|Du|dx+O_n.
  \end{aligned}
  \end{equation}

  From \eqref{defu1}, \eqref{esdef1}, \eqref{esdef2} and \eqref{esdef3} we deduce that
  \begin{equation}\label{defu1b}
  \begin{aligned}
    \mu\big(B(x_0,R_2)\big)\leq&\ \nu\big(\overline{B}(x_0,R_0)\big)
    +C\int_{B(x_0,R_0)} \big({h}^L+\nu^L+\varpi+(1+({\textsc{a}}^L)^\frac{1}{r})|D(u-v)|^p\big)^\frac{p-1}{p}({\textsc{a}}^L)^\frac{1}{pr}|D(u-v)|dx\\
    &+C({h}+\nu+\varpi)\big(\{R_0\leq|x-x_0|\leq R_1\}\big)\\
    &+C\int_{\{R_0<|x-x_0|<R_1\}}\big({h}^L+\nu^L+\varpi^L+(1+({\textsc{a}}^L)^\frac{1}{r})(|Dv|^p+|Du|^p)\big)^\frac{p-1}{p}({\textsc{a}}^L)^\frac{1}{pr}(|Du|+|Dv|)\,dx\\
    &+C({h}+\mu+\varrho)\big(\{R_1\leq|x-x_0|\leq R_2\}\big)\\
    &+C\int_{\{R_1<|x-x_0|<R_2\}}\big({h}^L+\mu^L+\varrho^L+(1+({\textsc{a}}^L)^\frac{1}{r})|Du|^p\big)^\frac{p-1}{p}({\textsc{a}}^L)^\frac{1}{pr}|Du|dx.
  \end{aligned}
  \end{equation}
Taking $R_0$ such that
\[
({h}+\nu+\varpi+\mu+\varrho)\big(\{|x-x_0|=R_0\}\big)=0,
\]
which holds true except for a countable set $E_{x_0}\subset \big(0, \dist(x_0,\partial\omega)\big)$, and making $R_1,R_2$ tend to $R_0$, from \eqref{defu1b} we deduce that
\begin{align*}
    \mu\big(B(x_0,R_0)\big)\leq&\ \nu\big(B(x_0,R_0)\big)\\
    &+C\int_{B(x_0,R_0)}\big({h}^L+\nu^L+\varpi^L+(1+({\textsc{a}}^L)^\frac{1}{r})|D(u-v)|^p\big)^\frac{p-1}{p}({\textsc{a}}^L)^\frac{1}{pr}|D(u-v)|dx,%\label{defu1c}
  \end{align*}
for any $R_0\in\big(0, \dist(x_0,\partial\omega)\big)\setminus E_{x_0}$ (observe that the right term in the integral is well defined as an element of $L^1(\omega)$). Therefore, the measures differentiation theorem shows \eqref{acopmunu}.

\cqfd
\par\bigskip\noindent
{\bf Proof of Lemma~\ref{lemafund2}.}
The proof is the same as the proof of Lemma~\ref{lemafund} choosing any point $x_0$ in $\Om$ rather than $\om$, extending the functions $u_n,v_n$ by $u$ in $\Om\setminus\om$, and then noting that the function $w_n$ defined by \eqref{wn} in $\Om$ is also equal to $u$ in $\Om\setminus\om$.
\cqfd
\par\bigskip\noindent
{\bf Proof of Corollary~\ref{clemafund}.}
     Assume that $1<p\leq N\!-\!1$. Applying Lemma \ref{lemafund} with $\omega=\omega_1$ (see also Remark~\ref{remlemafund} about the subsets of $\om$) we obtain
    \[\mu\leq\nu+C\big({h}^L+\nu^L+\varpi^L+(1+({\textsc{a}}^L)^\frac{1}{r})|D(u-v)^p| \big)^\frac{p-1}{p}({\textsc{a}}^L)^\frac{1}{pr}|D(u-v)|\quad\mbox{in }\om_1\cap\om_2.\]
    Analogously with $\omega=\omega_2$, we get
    \[\nu\leq\mu+C\big({h}^L+\mu^L+\varrho^L+(1+({\textsc{a}}^L)^\frac{1}{r})|D(u-v)^p| \big)^\frac{p-1}{p}({\textsc{a}}^L)^\frac{1}{pr}|D(u-v)|\quad\mbox{in }\omega_1\cap\omega_2.\]
    These two expressions prove the first estimate of \eqref{acomu-nu}.
    The proof of the second estimate is similar.
    \cqfd
\par\bs\noindent
{\bf Acknowledgement.}
The authors are grateful for support from the Spanish {\em Ministerio de Econom\'ia y Competitividad} through Project MTM2011-24457, and from the {\em Institut de Recherche Math\'ematique de Rennes}. The first author thanks the {\em Universidad de Sevilla} for hospitality during his stay April 18 - May 3 2016, and the second author thanks the {\em Institut de Math\'ematiques Appliqu\'ees de Rennes} for hospitality during his stay June~29 - July~10 2015.


\begin{thebibliography}{30}

\bibitem{Bal}{\sc J.M.~Ball}: ``Convexity conditions and existence theorems in nonlinear elasticity", {\em Arch. Rational Mech. Anal.}, {\bf 63} (1977), 337-403. 

\bibitem{BeBo}{\sc M. Bellieud \& G. Bouchitt\'e}: ``Homogenization of elliptic problems in a fiber reinforced structure. Nonlocal effects", {\em Ann. Scuola Norm. Sup. Pisa Cl. Sci.}, {\bf 26} (4) (1998), 407-436.

\bibitem{BeDe} {\sc A.~Beurling \& J.~Deny}: ``Espaces de Dirichlet", {\em Acta Matematica}, {\bf 99} (1958), 203-224.


\bibitem{Bra}{\sc A.~Braides}: {\em $\Gamma$-convergence for Beginners}, Oxford University Press, Oxford 2002, pp.~218.

\bibitem{BBC}{\sc A.~Braides, M.~Briane, \& J.~Casado D\'iaz}: ``Homogenization of non-uniformly bounded periodic diffusion energies in dimension two", {\em Nonlinearity}, {\bf 22} (2009), 1459-1480.


\bibitem{Bri}{\sc M.~Briane}: ``Nonlocal effects in two-dimensional conductivity", {\em Arch. Rat. Mech. Anal.}, {\bf 182} (2) (2006), 255-267.

\bibitem{BrCE1}{\sc M.~Briane \& M.~Camar-Eddine}: ``Homogenization of two-dimensional elasticity problems with very stiff coefficients", {\em J. Math. Pures Appl.}, {\bf 88} (2007), 483-505.

\bibitem{BrCa}{\sc M.~Briane \& J.~Casado D\'iaz}: ``Two-dimensional div-curl results. Application to the lack of nonlocal effects in homogenization", {\em Com. Part. Diff. Equ.}, {\bf 32} (2007), 935-969.

\bibitem{BrCa1}{\sc M.~Briane \& J.~Casado D\'iaz}: ``Asymptotic behavior of equicoercive diffusion energies in two dimension", {\em Calc. Var. Part. Diff. Equa.}, {\bfÊ29} (4) (2007), 455-479.

\bibitem{BrCa3}{\sc M.~Briane \& J.~Casado D\'iaz}~: ``Homogenization of convex functionals which are weakly coercive and not equibounded from above", {\em Ann. I.H.P. (C) Non Lin. Anal.}, {\bf 30} (4) (2013), 547-571.

\bibitem{BrCa4}{\sc M.~Briane \& J.~Casado D\'iaz}~: ``Homogenization of systems with equi-integrable coefficients", {\em ESAIM: COCV}, {\bf 20} (4) (2014), 1214-1223.

\bibitem{BrCa5}{\sc M. Briane \& J. Casado-Diaz:} ``A new div-curl result. Applications to the homogenization of elliptic systems and to the weak continuity of the Jacobian", to appear in {\em J. Diff. Equa.}

\bibitem{BCM}{\sc M.~Briane, J.~Casado D\'iaz \& F.~Murat}:  ``The div-curl lemma `trente ans apr\`es': an extension and an application to the G-convergence of unbounded monotone operators", {\em J. Math. Pures Appl.}, {\bf 91} (2009), 476-494.


\bibitem{BuDa}{\sc G. Buttazzo \& G. Dal Maso:} ``$\Ga$-limits of integral functionals", {\em J.~Analyse Math.}, {\bf 37} (1980), 145-185.

\bibitem{CESe1}{\sc M.~Camar-Eddine \& P.~Seppecher}: ``Closure of the set of diffusion functionals with respect to the Mosco-convergence", {\em Math. Models Methods Appl. Sci.}, {\bf 12} (8) (2002), 1153-1176.

\bibitem{CESe2}{\sc M.~Camar-Eddine \& P.~Seppecher}: ``Determination of the closure of the set of elasticity functionals", {\em Arch. Ration. Mech. Anal.}, {\bf 170} (3) (2003), 211-245.

\bibitem{CaSb}{\sc L.~Carbone \& C.~Sbordone:} ``Some properties of $\Ga$-limits of integral functionals",
{\em Ann. Mate. Pura Appl.}, {\bf 122} (1979), 1-60.

\bibitem{Cia}{\sc P.G.~Ciarlet:} {\em Mathematical Elasticity, Vol. I: Three-dimensional elasticity.} Studies in Mathematics and its Applications 20, North-Holland Publishing Co., Amsterdam 1988, pp.~451.

\bibitem{Cia2}{\sc P.G.~Ciarlet:} {\em Introduction \`a l'analyse num\'erique matricielle et \`a l'optimisation} (French) [Introduction to matrix numerical analysis and optimization], Math\'ematiques Appliqu\'ees pour la Ma\^{\i}trise [Applied Mathematics for the Master's Degree] Masson, Paris 1982, 279 pp.


\bibitem{Dac}{\sc B.~Dacorogna}: {\em Direct methods in the calculus of variations}, Applied Mathematical Sciences 78, Springer-Verlag, Berlin 1989, pp.~308.

\bibitem{Dal}{\sc G.~Dal Maso}: {\em An introduction to $\Ga$-convergence}, Birkha\"user, Boston 1993, pp.~341.

\bibitem{DeG1}{\sc E.~De Giorgi}: ``Sulla convergenza di alcune successioni di integrali del tipo dell'area", {\em Rend. Mat. Roma}, {\bf 8} (1975), 277-294.

\bibitem{DGFr}{\sc E. De Giorgi \& T. Franzoni:} ``Su un tipo di convergenza variazionale", {\em Rend. Acc. Naz. Lincei Roma}, {\bf 58} (6) (1975), 842-850.

\bibitem{FeKh}{\sc V.N.~Fenchenko \& E.Ya.~Khruslov}: ``Asymptotic behavior of solutions of differential equations with strongly oscillating matrix of coefficients which does not satisfy the condition of uniform boundedness", {\em Dokl. AN Ukr.SSR}, {\bf 4} (1981).

\bibitem{Fra}{\sc G.A.~Francfort:} ``Homogenisation of a class of fourth order equations with application to incompressible elasticity", {\em Proc. Roy. Soc. Edinburgh Sect. A}, {\bf 120} (1-2) (1992), 25-46.

\bibitem{Khr}{\sc E.Ya.~Khruslov}: ``Homogenized models of composite media", {\em Composite Media and Homogenization Theory}, ed. by G.~Dal Maso and  G.F.~Dell'Antonio, in {\em Progress in Nonlinear Differential Equations and Their Applications}, Birkha\"user 1991, 159-182.

\bibitem{Mos}{\sc U.~Mosco}: ``Composite media and asymptotic Dirichlet forms", {\em J.~Func. Anal.}, {\bf 123} (2) (1994), 368-421.

\bibitem{Mur}{\sc F.~Murat}: ``$H$-convergence", {\em S\'eminaire d'Analyse Fonctionnelle et Num\'erique}, 1977-78, Universit\'e d'Alger, multicopied, 34 pp. English translation~: {\sc F. Murat \& L. Tartar:} ``H-convergence", {\em Topics in the Mathematical Modelling of Composite Materials}, ed. by L.~Cherkaev \& R.V.~Kohn, Progress in Nonlinear Differential Equations and their Applications, {\bf 31}, Birka\"user, Boston 1998, 21-43.

\bibitem{APM}{\sc A. Pallares-Mart\'{\i}n:} ``High-Contrast homogenization of linear systems of partial differential equations", to appear in {\em Math. Meth. Appl. Sc.}

\bibitem{PiSe}{\sc C.~Pideri \& P.~Seppecher:} ``A second gradient material resulting from the homogenization of an heterogeneous linear elastic medium", {\em Continuum Mech. and Thermodyn.}, {\bf 9} (5) (1997), 241-257.

\bibitem{Rao}{\sc A. Raoult}: ``Nonpolyconvexity of the stored energy function of a Saint-Venant-Kirchhoff", {\em Material. Apl. Mat.}, {\bf 31} (6) (1986), 417-419.

\bibitem{San}{\sc E.~S\'anchez-Palencia:} {\em Nonhomogeneous Media and Vibration Theory}, Lecture Notes in Physics {\bf 127}, Springer-Verlag, Berlin-New York 1980, pp.~398.

\bibitem{Spa}{\sc S.~Spagnolo}: ``Sulla convergenza di soluzioni di equazioni paraboliche ed ellittiche", {\em Ann. Scuola Norm. Sup. Pisa Cl. Sci.}, {\bf 22} (3) (1968), 571-597.

\bibitem{Tem}{\sc R.~Temam:} {\em Probl\`emes math\'ematiques en plasticit\'e}, (French) [Mathematical problems in plasticity], M\'ethodes Math\'ematiques de l'Informatique [Mathematical Methods of Information Science] 12, Gauthier-Villars, Montrouge 1983, pp.~353.

\bibitem{Tar}{\sc L.~Tartar}: {\em The General Theory of Homogenization: A Personalized Introduction}, Lecture Notes of the Unione Matematica Italiana, Springer-Verlag, Berlin Heidelberg 2009, pp.~471.

\end{thebibliography}
\end{document}